%% file: main.tex
\documentclass{amsart}

\usepackage[english]{babel}
\usepackage{bbm}
\usepackage{hyperref}
\usepackage{color}
\usepackage{graphicx}
\usepackage{subfigure}
\usepackage[normalem]{ulem}

\newcommand{\N}{\ensuremath{\mathbb{N}}}
\newcommand{\R}{\ensuremath{\mathbb{R}}}
\newcommand{\C}{\ensuremath{\mathbb{C}}}

\newcommand{\E}{\ensuremath{\mathbb{E}}}
\renewcommand{\P}{\ensuremath{\mathbb{P}}}

\newcommand{\ind}[1]{\ensuremath{\mathbbm{1}_{\{#1\}}}}
\newcommand{\diff}{\mathop{}\mathopen{}\mathrm{d}}
\newcommand{\cal}[1]{\ensuremath{\mathcal{#1}}}

\newcommand{\eps}{\varepsilon}
\newcommand{\der}[2]{ \frac{\text{d} #1}{\text{d} #2} }  

\newcommand{\W}{\E(V)}
\newcommand{\eqdef}{\stackrel{\text{def.}}{=}}

\newcommand\croc[1]{\left\langle #1\right\rangle}
\newcommand\Lam[2]{\croc{\Lambda(#1),\, #2}}
\newcommand\Delt[2]{\croc{\Delta(#1),\, #2}}

\definecolor{orange}{rgb}{1.00,0.50,0.0}

\graphicspath{{./Figures/}}

\definecolor{orange}{rgb}{1.00,0.50,0.0}

\def\etal{et al.}

\newtheorem{proposition}{Proposition}
\newtheorem{lemma}{Lemma}
\newtheorem{theorem}{Theorem}

\title[Random Integrate-and-Fire neurons]{On the dynamics of random neuronal networks}
\date{\today}
\author[Ph. Robert]{Philippe Robert}
\email{Philippe.Robert@inria.fr}
\address[Ph. Robert]{INRIA Paris---Rocquencourt, Domaine de Voluceau, 78153 Le Chesnay, France. }
\urladdr{\href{http://team.inria.fr/rap/robert}{http://team.inria.fr/rap/robert}}
\author[J. Touboul]{Jonathan Touboul}
\email{jonathan.touboul@college-de-france.fr}
\address[J. Touboul]{Mathematical Neuroscience Team, CIRB - Collège de France and INRIA Paris-Rocquencourt\\11, Place Marcelin Berthelot
75005 Paris, FRANCE }
\urladdr{\href{http://mathematical-neuroscience.net/team/jonathan/}{http://mathematical-neuroscience.net/team/jonathan/}}
\begin{document}

\begin{abstract}
We study the mean-field limit and stationary distributions of a pulse-coupled network modeling the dynamics of a large neuronal assemblies. Our model takes into account explicitly the intrinsic randomness of firing times, contrasting with the classical integrate-and-fire model.  The ergodicity properties of the Markov process associated to finite networks are investigated. We derive the limit in distribution of the sample path of the state of a neuron of the network  when its size gets large.   The invariant distributions  of this limiting  stochastic process are analyzed as well as their stability properties. We show that the system undergoes transitions as a function of the averaged connectivity parameter, and can support trivial states (where the network activity dies out, which is also the unique stationary state of finite networks in some cases) and self-sustained activity when connectivity level is sufficiently large, both being possibly stable. 
\end{abstract}

\maketitle

\begin{center} {\bf Work in progress - Preliminary version of \today}
\end{center}
\bigskip

\hrule

\vspace{-3mm}

\tableofcontents

\vspace{-1cm}

\hrule

\input Intro

\input Model
\input Finite
\input McKeanVlasov

\input Mean-Field
\input Invariant
\input Appendix
\input AppendixSimu

\bibliographystyle{amsplain}
\bibliography{ref}

\end{document}

%% file: Intro.tex
\section{Introduction}
\subsection{Leaky Integrate and Fire Neuron Models}
We investigate a model of neuronal network in which the state of a neuron $i$ at time $t$ is given by a jump process describing the membrane potential of a cell. Transitions occur either when a cell receives an action potential from another cell in the network, or when the cell fires a spike and is reset to its resting potential. This event occurs randomly at a state-dependent rate. The inspiration for the development of this model is the celebrated leaky integrate-and-fire neuron (see biological motivation below in section~\ref{sec:Motiv}). The original integrate-and-fire model has been studied under the assumptions that the membrane potential of each cell is noisy, and that spikes are fired as soon as a fixed threshold $V_F$ is  exceeded. It can be described as an hybrid ODE with jumps: 
\begin{itemize}
\item   if ${ X_i(t-)<V_F}$\\
\[
\diff X_i(t)=-X_i(t)\diff t+\sum_{j\not=i}  w_{ji}(t) \diff D_j(t)+ \diff I_i^e(t),
\]
\item  if ${ X_i(t-)=V_F}$\\
 $\diff D_i(t) =1,$ $X_i(t)=V_R$,
\end{itemize}
with $f(t-)$ the left-limit of $f(\cdot)$ at $t$ and  $D_j(t)$ the number of spikes of cell $j$ up to time $t$ and $w_{ji}(t)$ is the input of cell $j$ sent to cell $i$ when it fires at time $t$. The quantity $I_i^e(t)$ is the external input to cell $i$.

Caceres, Carrillo and Perthame in~\cite{caceres-carrillo:11} have demonstrated the interesting property that in the mean-field limit, the solutions blow-up in finite time. While this may attributed to the instantaneity of spikes firings and their immediate effects on the firing of other cells, it remains a non-trivial issue since this blow-up also occurs when considering propagation delays~\cite{caceres2014beyond}, which prevents from avalanche phenomena and instantaneous self-excitation. 

It should be noted that in most of works of the literature the stochastic component of this dynamical system is generally represented by the external input $(I_i^e(t))$, see Section~\ref{sec:bibli}. It is classically considered to be a Gaussian process describing the effect of the neurons belonging to other populations. In these models the firing mechanism of a given neuron is a deterministic function of the input received. Experimental data suggest nevertheless that the firing process has a significant random component. A number of theoretician and experimentalists have discussed this issue~\cite{chichilnisky2001simple,pillow2005prediction}, and we briefly summarize their findings in section~\ref{sec:Motiv}. All in all, a more accurate description incorporates the intrinsic spike time variability by considering that neurons fire according to a non-homogeneous Poisson process, with a rate being an increasing function of their voltage.

The stochastic model investigated in this paper describes the occurrences of spikes as follows. A  cell with membrane potential $x$ fires at rate $b(x)$ where $x\mapsto b(x)$ is a non-decreasing function. Note that the threshold model corresponds to the case $b(x){=}0$ if $x{<}V_F$ and $b(x){=}+\infty$ otherwise. In this context, the analogue of the above differential equations  are given by 
\begin{equation}\label{eqintro}
\diff X_i(t)={-}X_i(t)\diff t{+}\sum_{j\not=i}  W_{ji}(t-) {\cal N}^j_{b(X_j(t-))}(\diff t)\\ -X_i(t-){\cal N}^i_{b(X_i(t-))}(\diff t)
\end{equation}
where 
\begin{itemize}
\item  ${\cal N}^j_y(\diff t)$ denotes a Poisson point process with rate $y\geq 0$,
\item For $s\geq 0$, $(W_{ij}(s))$ is an i.i.d. sequence of non-negative integrable random variables,
\end{itemize}
and $f(t-)$ denotes the left limit of the function $f$ at $t$.
In particular, $V_R$ is set to be $0$, after a spike the state of a neuron  is $0$. A more formal formulation of these equations is given in Section~\ref{sec:model}. Note that intrinsic randomness occurs in this model through the firing mechanism with the Poisson processes but also through the fluctuations of synaptic efficacy $W_{ji}(t)$ representing the effect of a spike of neuron $j$ on the membrane potential of neuron $i$ at time $t$. 

The quantity $b(0)$ is the rate at which a neuron fires when its membrane potential is at rest, i.e. when it is not driven by spikes from other neurons. It can be interpreted as the firing rate of the cell due to external noise. In this way the firing rate can be decomposed as $b(x){=}b(0){+}[b(x){-}b(0)]$, the sum of external and intrinsic firing rates. In other words, our network incorporates input from cells not belonging to the network considered as a Poisson Process with constant rate: this is the classical interpretation behind the Gaussian input authors generally consider.

\subsection{Brief overview of the literature in neural mean-field dynamics}\label{sec:bibli}
As we have already remarked, networks of integrate and fire neurons have been mainly studied mathematically with the assumption that the randomness is essentially contained in the external input $(I^e_i(t))$ and that it is a continuous process driven by a Brownian motion. In this setting the classical tools of stochastic calculus, It\^o Formula in particular for probabilistic approaches, can be used to investigate the behaviour of these networks. 

There are several works using PDEs, see Caceres \etal~\cite{caceres-carrillo:11,caceres2014beyond}, Pakdaman \etal~\cite{pakdaman2010dynamics,pakdaman2012adaptation,pakdaman2013relaxation}. Recently, these were complemented by probabilistic approaches Delarue \etal~\cite{delarue:hal-01001716,tanre:12} and Inglis and Talay~\cite{inglis2014mean} with Gaussian processes. At the other end of the spectrum of mathematical models for neurosciences are the detailed conductance-based models. These systems represent accurately the different ionic exchanges at play during the firing of an action potential. This has been the subject of a number of recent works in the domain, see~\cite{touboulneuralfields:14,touboulSingular:14} in the context of spatially extended networks. Although these models lead to complex equations, some simple models for the firing rate of neurons, are tractable. For these networks, it should be noted that there is an interesting behavior as a function of noise levels in a mean-field context, see~\cite{touboulNeuralFieldsDynamics:12}.

To the best of our knowledge mathematical models of neural networks using Poisson processes to describe the occurrences of firing events are more rare compared to models with Brownian random components. The main reason is, probably, the fact that the corresponding stochastic calculus is more awkward to deal with. But there has been a renewed interest recently. De Masi \etal~\cite{demasi-presutti:14} investigated a model with firing rates depending on the state of the cell but without the leaking term of the integrate-and-fire model. They proved an hydrodynamic limit behavior by using coupling methods on a discrete-time version of the process. The same model has been investigated by Fournier and L{\"o}cherbach~\cite{FL2014} where a quite general mean-field convergence result is proved with an explicit upper bound on the rate of convergence to the mean-field model. In this paper~\cite{FL2014}, appeared independently at the same time as the first version of the present paper, the analysis of Section~3 and~4 parallel, with different methods, Sections~\ref{McKeanSec} and \ref{sec:Convergence} below.  

\subsection{Summary of the results}
The main goal of the present paper is to investigate the stability properties of these networks whose firing mechanisms are driven by Poisson processes. The state of the network is described in terms of the solutions of Equation~\eqref{eqintro}. We first consider finite-sized networks and then the limit of the system as the network size tends to infinity; in both case, we analyze of the corresponding invariant distributions. 

\medskip
\noindent
{\bf Finite Networks}\\
It is shown that if there is an external noise, i.e. $b(0){>}0$, then there exists a unique non-trivial invariant distribution for the Markov process $(X_i(t))$. An explicit representation of the invariant distribution is given when $b(\cdot)$ is a constant function.  When $b(0){=}0$ it is shown that, almost surely no spike occurs after some time. The network \emph{dies} with probability $1$ in this case, and the Dirac mass at $0$ is the unique, trivial, invariant distribution. 

The Harris recurrence properties of the associated Markov process are proved in Section~\ref{sec:finite}. They rely on a regeneration argument and the invariant distribution is obtained via a backward coupling argument.

\medskip
\noindent
{\bf Large Networks}\\
For large size networks, the output of cell $j$ to cell $i$ when cell $j$ fires at time $t$, the variable $(W_{ji}(t))$, is taken as $(V_{ji}(t)/N)$ where $(V_{ji}(t))$ are i.i.d. integrable random variable so that Equation~\eqref{eqintro} of the evolution of $(X_i^N(t))$ becomes
\begin{equation}\label{eqintro2}
\diff X^N_i(t)={-}X^N_i(t)\diff t{+}\frac{1}{N}\sum_{j\not=i}  V_{ji}(t-) {\cal N}^j_{b(X^N_j(t-))}(\diff t)\\ -X^N_i(t-){\cal N}^i_{b(X^N_i(t-))}(\diff t)
\end{equation}
It is shown that provided that $b(\cdot)$ satisfies some growth condition and other technical conditions, a mean-field result holds. In particular,  the distribution of the sample path of the state of a cell of the network converges in distribution to the solution of the non-linear stochastic differential equation, the associated McKean-Vlasov process $(Z(t))$
\begin{equation}\label{eqintro3}
\diff Z(t)=-Z(t)\diff t +E(V)\E(b(Z(t)))\diff t -Z(t-) {\cal N}_{b(Z(t-))}(\diff t). 
\end{equation}
One of the problems here is to cope with the fact that the function $(b(x))$ is not necessarily Lipschitz like $b(x)=\lambda x^\alpha$ with $\alpha>1$ which will be considered. This is generally a source of trouble to get mean-field results. See Sznitman~\cite{Sznitman} page~176 and Scheutzow~\cite{Scheutzow} for example. As it will be seen, a control of the sample paths of the firing events is the main technical ingredient to overcome this problem.   Another difficulty is to deal with the discontinuities due to Poisson events and the associated It\^o calculus. 

\medskip
\noindent
{\bf Invariant States}\\
The last part of the paper is devoted to the properties of the invariant distributions of Equation~\eqref{eqintro3}.  When $b(0){=}0$ the Dirac mass at $0$, $\delta_0$,  is always a possible invariant distribution, its stability properties are investigated. 

In the linear case, $b(x){=}\lambda x$, it is shown that if $\lambda\E(V){<}1$ then $\delta_0$ is the unique invariant distribution but if $\lambda\E(V){>}1$ there exists a non-trivial invariant distribution for $(Z(t))$.
The existence of non-trivial distribution suggests a quasi-stationary phenomenon. It is known that, in this case,  the state of the network of finite size is dying with probability $1$ absorbed at $0$. Nevertheless if $\lambda\E(V){>}1$,  starting from a large initial state, the time before absorption at $0$ is quite large. Before absorption, the state of the network can be, for a while, in a local equilibrium. 

More striking, it is  shown that if $b(x){=}\lambda x^\alpha$ with $\alpha{>}1$, then there exists some $\rho_c{>}0$ such that if $\lambda\E(V){>}\rho_c$ then, in addition to the trivial invariant distribution, there exist at least {\em two} non-trivial invariant distributions for $(Z(t))$. Simulations presented at the end seem to suggest that one of them is stable but not the other one and that $\delta_0$ has a non-empty basin of attraction. 

\bigskip

\noindent
{\bf Outline of the paper}\\
The paper is organized as follows. In section~\ref{sec:model} we present our model and the main results of the paper, before addressing in section~\ref{sec:finite}  the properties of finite networks. The following sections are devoted to the analysis of the limit as the system size diverges: we analyze the properties of the McKean-Vlasov limit in section~\ref{McKeanSec}, and show the mean-field convergence of the network towards this limit in section~\ref{sec:Convergence}. We conclude with the analysis of invariant distributions of the McKean-Vlasov equation in section~\ref{sec:Convergence}.

\medskip
\noindent
{\bf Acknowledgments}\\
The authors  are grateful to Nicolas Fournier for his remarks and for pointing out a mistake in the proof of Theorem~\ref{ThMF} in the first version of the paper. 

\section{Biological background}\label{sec:Motiv}
Neurons are intrinsically noisy electrically excitable cells, that transmit information through stereotyped electrical impulses, called action potentials, or spikes. Spikes are transmitted through synapses to all connected neurons, which have the effect of either increasing (excitatory) or decreasing (inhibitory) their membrane potential. In the cortex, neurons tends to form large populations of statistically identical cells, receiving the same input and strongly interconnected. These cortical areas (or cortical columns) are of the order of a few millimeters and contain hundreds to hundred of thousands of neurons. In such structures, neurons fire spikes when the difference of electrical potential between the intra- and extra-cellular domains is sufficiently large (hyperpolarization), therefore occurring in response to a sufficient amount of excitatory spikes received by the neuron. 

The spiking nature of the neuronal activity motivated the introduction of a simple heuristic model of single cell, the \emph{integrate-and-fire} neuron, which makes the assumption that the membrane potential linearly integrates the inputs it receives and fires a spike as soon as a fixed voltage threshold is reached. This model, introduced in the beginning of last century by Louis Lapique~\cite{lapicque:07}, has given rise to an extensive literature (see e.g. Burkitt~\cite{burkitt:06,burkitt:06a} for reviews). The introduction of this model was instrumental in the advances made in the mathematical and computational understanding neurons' activity. Though overly simplified, this model has allowed to take advantage of the stereotyped, fast nature of spikes to focus on the important problem of spike timings. Moreover, this model, for its relative simplicity, has allowed interesting mathematical developments and yielded better understanding of isolated cells, see Brunel~\cite{brunel:00}.

In biological conditions, it has been observed that neurons display a noisy activity. The first model integrating randomness is due to Gerstein and Mandelbrot~\cite{gerstein-mandelbrot:64} who incorporated to the model random spike arrivals as a random walk. This model did not considered the intrinsic nature of the nerve cells firing times, but rather considered that this randomness was due to incoming spikes from cells outside of the network considered, i.e. disregarded the fact that isolated cells actually display an intrinsically noisy activity, i.e. fire irregularly even when disconnected from their network. This seminal model, interesting in many regards, lead to many developments. In particular, diffusion limits of the incoming spike strain using stochastic differential equations, see Knight~\cite{knight:72} and Stein~\cite{stein:65} lead to the introduction of the celebrated \emph{leaky integrate-and-fire model}. This type of model has been paramount in the study of noisy integrate-and-fire since then and lead to develop new technique to analyze their singular behavior.  See e.g. Brunel~\cite{brunel:00,brunel:00b}, Delarue \etal~\cite{tanre:12} and Caceres \etal~\cite{caceres-carrillo:11} to cite a few. 

At the level of one cell, classical models assume random arrivals of spikes. In a diffusion limit, it is represented as an external additive Brownian noise.   In the context of a large set of cells, randomness of spike times due to intrinsic variability is nevertheless an important phenomenon which cannot be ignored, see e.g. Rolls and Deco~\cite{rolls-deco:10}. Another drawback of the classical noisy integrate-and-fire neuron is the presence of a fixed threshold. This feature, abstraction of the actual dynamics of neurons, induces a number of artifacts, among which the so-called avalanche effects. It occurs in excitatory networks, corresponding to an explosion of the spike rate: all neurons fire instantaneously at the same time inducing again all neurons to fire again. And these phenomena do occur, generically, in the mathematical models of integrate-and-fire networks, as shown in Caceres \etal~\cite{caceres-carrillo:11}. All this context points to the fact that the integrate-and-fire neuron model, interesting in order to understand single isolated cells in a non-noisy context, may not be the best model when considering the behavior of very large networks of noisy neurons as appearing in physiological conditions in the brain.

Accounting for the intrinsic randomness of firing times has been a longstanding issue in computational neurosciences. The fact that nerve cells integrate the input received and fire spikes at random times with an intensity depending on their membrane potential. Instead of continuous Gaussian processes, a natural model is to consider that the membrane potential is driven by inhomogeneous Poisson jumps whose intensity is a function of the voltage of the cells. This model has been introduced one century later than the integrate-and-fire model in Chichilnisky~\cite{chichilnisky2001simple} and is called linear-nonlinear Poisson model. This model seems particularly well suited in order to represent the neuronal firing, and displayed a good fit with experimental data, allowing precise prediction of spike trains Pillow \etal~\cite{pillow2005prediction,pillow2008spatio}. The present study use this model as the building block of our networks. 

The model we shall study in the present paper has the interest of conserving the discrete nature of spikes, allowing to characterize the statistics of spike trains in the limit of large networks. Moreover, taking into account intrinsic noise, beyond its biological interest, ensures well-posedness of the system, allowing to describe the limit without resorting to additional processes. Eventually, an interesting property is that the dynamics of the limit equation can be partially characterized, and stationary solutions can be described. As it will be seen, depending on the sharpness of the spiking intensity function as well as the average coupling strength, self-sustained spontaneous activity may arises in infinite networks. This original phenomenon opens interesting questions on these models.

%% file: Model.tex
\section{Stochastic Model}\label{sec:model}
We consider a network composed of $N$ neurons, whose state is described by a scalar variable representing its membrane potential $X_i^N(t)$ (i.e., the difference of electric potential between the intra- and extra-cellular domains). This quantity decays exponentially fast towards zero in the absence of input due to the leak currents and ions flowing across the cellular membrane, leading the voltage to its equilibrium value, assumed here to be $0$. The timescale of this process is our time unit, i.e., for $1\leq i\leq N$, if the neuron does not spike and does not receive any spike in the interval of time $[T,T']$, the membrane potential of neuron $i$ satisfies the ordinary differential equation:
\[\der{X_i^N(t)}{t}=-X_i^N(t).\]

Neurons fire at random times, according to a voltage-dependent Poisson process with rate $b(x)$ where $x$ denotes the voltage of the neuron. After spiking, the neuron's voltage is instantaneously reset to its rest potential $X_i(t)=0$, and the voltage of neurons $j\not=i$ are instantaneously updated: their voltage is added the synaptic coefficient $W_{ij}$, which are considered to be i.i.d. random variables with law $F$:
\[X_j(t)=X_j(t-)+W_{ij}.\] 
For simplicity, it is assumed that the random variables  $(W_{ij},j\not=i)$ are positive i.i.d. and integrable and that their distribution does not depend on $i$. In this way, it is easily seen that the process $(X(t))$ is Markov process. If $x=(x_i)\in\R_+^N$, $\|x\|$ denotes the $l^1$ norm:
$\|x\|= \vert x_1\vert+\cdots+\vert x_N\vert$.
\subsection*{Evolution Equations}
An equivalent description of $(X(t))$ can be provided in terms of the solution of the following Stochastic Differential Equation (SDE),
\begin{multline}\label{SDE}
\diff X_i(t)=-X_i(t)\diff t +\sum_{j\not=i} \int_{\R_+^2} z_i\ind{0\leq u\leq b(X_j(t-))} {\cal N}_j(\diff u,\diff z,\diff t)
\\-X_i^N(t-)\int_{\R_+^2} \ind{0\leq u\leq b(X_j(t-))} {\cal N}_j(\diff u,\diff z,\diff t),
\end{multline}
where $({\cal N}_j)$ are independent Poisson processes,  for $1\leq i\leq N$,  ${\cal N}_i$ has the intensity measure  given by $ \diff u \otimes {\mathbbm{W}^i(\diff z)} \otimes \diff t $
where \[{\mathbbm{W}^i(\diff z)= \otimes_{j=1}^{i-1} W(\diff z_j)\otimes\delta_0(\diff z_i)\otimes_{j=i+1}^N W(\diff z_j)}\]
is the measure corresponding to the result of the emission of a spike by neuron $i$ on the voltage of all neurons. In the latter expression, $\delta_0$ is the Dirac distribution at $0$ and $W(\diff x)$ is the common distribution of the random variable $(W_{ij}, j\not=i)$ on $\R_+$ associated with the amount of excitation received by a neuron after a spike of another neuron. 

Equivalently, it can be written as
\begin{multline}\label{SDE1}
 X_i(t)=-\int_0^t X_i(s)\,\diff s +\E(W_1) \sum_{j\not=i}\int_0^t  b(X_j(s))\,\diff s
\\-\int_0^t X_i^N(s)b(X_j(s))\diff s+M_i(t),
\end{multline}
where 
\begin{multline}\label{SDE-mart}
M_i(t){=}\sum_{j\not=i}\int_{s=0}^t \int_{u=0}^{b(X_j(s-))} 
\left[\int_{z=0}^{+\infty} z_i {\cal N}_j(\diff u,\diff z,\diff s){-}\E(W_1)\diff u \diff s\right]
\\-\int_{s=0}^t   X_i(s-) \int_{u=0}^{b(X_i(s-))}
\left[\int_{z=0}^{+\infty}  {\cal N}_i(\diff u,\diff z,\diff s){-}\diff u \diff s\right] ,
\end{multline}
is the associated local martingale. See Proposition~\ref{PoisMart} of the appendix and Rogers and Williams~\cite{rogers} for example.

\subsection*{Extinction properties of the network}
From the biological viewpoint, it is natural to assume that $x\mapsto b(x)$  is a non-decreasing positive function since the higher the potential, the more likely a spike will occur. If the initial state of a neuron is $x$ and if no spike occurs in the network (no neuron fires) in the time interval $[0,t]$, its state at time $t$ is equal to $x\exp(-t)$, and in particular its instantaneous firing rate at this time, $b(x\exp(-t))$ decreases to $b(0)$ if $t$ diverges. If this later quantity is $0$, it may happen that  the neuron will not spike with positive probability. In this case, if the components  $x_i(0)$, $1\leq i\leq  N$  of the  initial value of the state of the network are too small, there would be an event of positive probability  for which no spike occurs at all.  The following lemma provides a sufficient condition on the behavior of the map $b(x)$ at $0$ under which extinction of the network activity does not occur. Theorem~\ref{Th-dead} below completes this result. 
\begin{lemma}\label{Prop-death}
If the condition
\begin{equation}\label{Cond}
 \int_{[0,1]} \frac{b(s)}{s}\,\diff s=+\infty
\end{equation}
holds, then a node with a non-zero initial value spikes with probability $1$. 
\end{lemma}
\begin{proof}
The function $x\mapsto b(x)$ being non-decreasing,  Relation~\eqref{Cond} also holds when $[0,1]$ is replaced by $[0,a]$ with $a>0$. Let $1\leq i\leq N$ and $X_i(0)=x>0$ then, if $\tau_i$ denote the instant (possibly infinite) of the first spike of neuron $i$, one has 
\begin{multline*}
\P(\tau_i>t)=\E\left(\exp\left(-\int_0^t b(X_i(s)\,\diff s\right)\right)\\\leq 
\exp\left(-\int_0^t b\left(xe^{-s}\right)\,\diff s\right)=\exp\left(-\int_{xe^{-t}}^x \frac{b(s)}{s}\,\diff s\right)
\end{multline*}
since the relation  $X_i(t)\geq x\exp(-t)$ holds for $t\leq \tau$ (the other neurons may only increase the state of neuron $i$) and that the function $x\mapsto b(x)$ is non-decreasing. By letting $t$ go to infinity, one gets that $\P(\tau_i=+\infty)=0$ by Condition~\eqref{Cond}. 
\end{proof}
Condition~\eqref{Cond} together with the monotonicity  and  a convenient regularity property  imply in fact that $b(0)$ is positive. The quantity $b(0)$ is the firing rate of a neuron with a flat potential, it can be see as a representation of the  external noise. 

We now investigate  the stability of the Markov process $(X_i(t))$. In the absence of spikes, each of the components decreases exponentially to $0$ and  it is reset to $0$ when the corresponding node fires. The ergodicity property seems to be quite likely provided that it is proved that the nodes do not fire too quickly as in the PDE description of  Caceres \etal~\cite{caceres-carrillo:11}. The analysis of these properties for finite-size networks are now investigated in order to ensure that these properties hold. 

%% file: Finite.tex
\section{Finite Networks}\label{sec:finite}
In this section the number of neurons will be kept fixed, so we drop the upper index $N$ throughout the section for simplicity of the notations.

\subsection{Recurrence, Ergodicity and Invariant Measures}
We start with a technical result related to an estimation of the mean return time in a specific compact set. 
\begin{proposition}\label{Prop-Foster}
There exists $C_0$ such that if
\[
T_0=\inf\left\{u>0: X(u)\in[0,C_0]^N\right\}
\]
then, for  $X(0)=x=(x_i)\not\in[0,C_0]^N$, 
\begin{equation}\label{Foster}
\E_x(T_0)\leq \|x\|=x_1+\cdots+x_N.
\end{equation}
\end{proposition}
\begin{proof}
Let
\[
F=\left\{x\in\R_+^N: \sum_{i=1}^N x_i[1+b(x_i)]\leq (N{-}1)\E(W_1)\sum_{i=1}^Nb(x_i) +1 \right\}
\]
due to the monotonicity property of the function $x\mapsto b(x)$,  $F$ is a compact subset of $[0,C_0]^N$, with $C_0=N^2\E(W_1)+1$. 
Denote
\[
T_F=\inf\left\{u>0: X(u)\in F\right\},
\]
then clearly $T_0\leq T_F$. 

If $X(0)=x\not\in F$ and $t\geq 0$, define $S(t){=}\|X(t)\|{=} X_1(t)+X_2(t)+\cdots+X_N(t)$, then Relation~\eqref{SDE1} gives the identity
\begin{multline*}
 S(t) = S(0)+ \int_0^t\left[ (N{-}1)\E(W_1) \sum_{i=1}^N  b(X_j(u))
\right. \\ \left. -\sum_{i=1}^N X_i(u)(1+b(X_i(u)))\right]\diff u +\sum_{i=1} M_i(t),
\end{multline*}
where $(M_i(t))$ are the local martingales defined by Equation~\eqref{SDE-mart}.

Assume that $X(0)\not\in F$, since $T_F$ is a stopping time, one has
\begin{align*}
0\leq \E(S(T_F&\wedge t))= S(0)\\
 &{+}\E\left(\int_0^{t\wedge T_F}\left[ (N{-}1)\E(W_1) \sum_{i=1}^N b(X_j(u))
{-}\sum_{i=1}^N X_i(u)[1{+}b(X_i(u))]\right]\diff u \right)\\
&\leq S(0) -\E(t\wedge T_F). 
\end{align*}
One gets $\E(t\wedge T_F)\leq {S(0)}$ and consequently the desired relation. 
\end{proof}

To state the stability properties of the Markov process $(X(t))$, the framework of Harris Markov processes now is used. See Nummelin~\cite{Nummelin} and Asmussen~\cite{Asmussen} for a general introduction. 
\begin{theorem}[Stability] \label{Th-dead}\ 
\begin{enumerate}
\item If the condition
\begin{equation}\label{M}
\int_0^{1} \frac{b(s)}{s}\,\diff s<+\infty, 
\end{equation}
holds then, almost surely, no spike occurs after some finite time, in particular
\[
\lim_{t\to+\infty} (X_i(t),1\leq i\leq N)=0
\]
and the Dirac mass at $0$ is the unique invariant distribution. 
\item If  $b(0)>0$ and if there exists $K>0$ such that $W_1\leq K$ a.s.,  then  the Markov process $(X_i(t), 1\leq i\leq N)$ is Harris ergodic. 
\end{enumerate}
\end{theorem}
\begin{proof}
Assume  that Condition~\eqref{M} holds. The notations of Proposition~\ref{Prop-Foster} are used. 
Let $X(0)=x\in\R_+^N$, $x\not=0$. If $X(0)=x\in [0,C_0]^N$, in the proof of Proposition~\ref{Prop-death} it has been seen that if $\tau_i$ is the first time neuron~$i$ spikes, $1\leq i\leq N$,  then in absence of spikes of the other nodes,
\[
\P(\tau_i=+\infty) =\exp\left(-\int_{0}^{x_i} \frac{b(s)}{s}\,\diff s\right)\geq\eta\eqdef  \exp\left(-\int_{0}^{C_0} \frac{b(s)}{s}\,\diff s\right)>0.
\]
 Consequently if $X(0)\in [0,C_0]^N$ there is a positive probability lower-bounded by $\eta^N$ that none of the nodes spike. From Proposition~\ref{Prop-Foster} one gets that if $(X(t))$ leaves $F$ then it returns with probability $1$, hence almost surely the process $(X(t))$ will stop having upward jumps after some time, (1) is proved. 

\medskip
Now it is assumed that $b(0)>0$.  The strategy of the proof is as follows,  in absence of spikes of other nodes, the duration of time for the next spike  of a node with  value $y>0$ can be decomposed as the minimum of two random variables with one of them not depending of $y$. This property provides a way of having a regeneration mechanism (i.e. forgetting the value $y$).  If all nodes proceed along the same line, then the initial value of the Markov process is forgotten after some time on some event of positive probability, which gives the key regenerative structure of a Harris Markov process.

If $X(0)=x=(x_i)$, one denotes by $E_1^{i}$ an exponential random variable with parameter $b(0)$ and $E_1^{x_i,i}$ a random distribution such that
\[
\P(E_1^{x_i,i}\geq t)=\exp\left(-\int_0^t \left[b\left(x_ie^{-u}\right)-b(0)\right]\,\diff u\right).
\]
Hence, starting from the initial state $x$, if no other spike occurs before, the first instant when node $i$ spikes has the same distribution as $E_1^{i}\wedge E_1^{x_i,i}$.  The variable $E_1$ is independent of the variables $(E^{x,i}_1, x>0)$ which can be chosen so that $E^{x,i}_1\leq E^{y,i}_1$ for $x\geq y$. 

For $n\geq 1$, one denotes by $Y_n=(y_{i,n})$ the state of the Markov process $(X(t))$ just after the $n${\em th} jump/spike. The sequence $(Y_n)$ is the embedded Markov chain, for $n\geq 1$ $Y_n$ is the state of $(X(t))$ at the instant of the $n$th jump.  Let $f$ be some non-negative Borelian function on $\R_+^N$. If $X(0)=x\in F$,  then, since $x_i\leq C_0$, 
\[
\E_x(f(Y_1))\geq \E\left(f\left(Y_1\right)\ind{E_1^1\leq \min_i(E_1^{i}\wedge E_1^{x_i,i})}\right) \geq  \E\left(f\left(Y_1\right) \mathbbm{1}_{A_1}\right),
\]
where 
\[
A_1=\left\{E_1^1\leq \min_{1\leq i\leq N}E_1^{i}\wedge E_1^{C_0,i}\right\}
\]
and, on $A_1$, 
\[
Y_1=(y_{1,1},y_{2,1},\ldots, y_{N,1})=(0,x_2e^{-E_1^1}{+}W_{1,2},\ldots, x_Ne^{-E_1^1}{+}W_{1,N}).
\] 
This inequality is associated to the event that node~1 spikes first on its ``$b(0)$-component'' $E_1^1$. As a result the lower bound of the above relation does not depend on $x_1$ anymore.  We proceed in the same way for the second step,  with node~2 spiking this time, and since $y_{i,1}\leq C_0+K$ for $1\leq i\leq N$, 
\begin{equation}\label{M1}
\E_x(f(Y_2))\geq \E\left(f\left(Y_2\right)\mathbbm{1}_{A_1\cap A_2}\right),
\end{equation}
with 
\[
A_2=\left\{E_2^2\leq \min_{1\leq i\leq N} E_2^{i}\wedge E_2^{C_0+K,i}\right\}
\]
and 
\[
Y_2=(W_{2,1},0,y_{3,1}e^{-E_2^2}{+}W_{2,3}, \ldots, y_{N,1}e^{-E_2^2}{+}W_{2,N}).
\] 
on the event $A_1\cap A_2$. This time the lower bound~\eqref{M1} does not depend on $x_1$ and $x_2$.  We can proceed recursively, and finally get the relation
\[
\E_x\left(f(Y_N))\geq \E(f(Z)\mathbbm{1}_{\cal A}\right),\quad \forall x\in F,
\]
where the random variable $Z$ and the set ${\cal A}$  do not depend on $x\in F$ and that $P({\cal A})>0$. Consequently $F$ is {\em a regeneration set} of the Markov chain $(Y_n)$,  see Asmussen~\cite[page~198]{Asmussen} for example. 
The Harris property of the Markov chain has been established. To prove the ergodicity, it is enough to prove that if
\[
T_0^+=\inf\left\{u>0: X(u)\in [0,C_0]^N \text{ and }\exists v\leq u, X(v)\not\in [0,C_0]^N \right\},
\]
is the first return time to $F$ after an exit, then 
\[
\sup_{x\in F}\E_x(T_0^+)<+\infty.
\]
See, for example, Asmussen~\cite[Theorem~3.2, page~200]{Asmussen} and Robert~\cite[Proposition~8.12, page~221]{Robert}.
For $x\in F$, the first time $(X(t))$ is in $y=(y_i)$ outside $F$, necessarily $y_i\leq C_0+K$, consequently, by the strong Markov property, 
\[
\E_x(T_0^+) \leq \sup_{\substack{y=(y_i)\not\in F,\\  \max_i y_i\leq C_0+K}} E_y(T_0)\leq  C_0+K,
\]
by Proposition~\ref{Prop-Foster}. The theorem is proved. 
\end{proof}

\subsection{The State-Independent Process}
The case when the firing rate function $b$ is constant is investigated. In this setting, the neurons spike independently of their state. This is one of the very rare cases where one can get some substantial information on the distribution of the equilibrium of the network. 
\begin{proposition} If the firing rate is constant  and equal to $\lambda>0$, then the invariant distribution of the Markov process $(X(t))$ is the law of the vector $(X_1,\ldots,X_N)$ with
\begin{equation}\label{mi-const}
X_i=\sum_{j\not=i} \sum_{k\geq 1} W^j_{ik}e^{-t_{jk}} \ind{t_{jk}\leq t_{i1}},\quad 1\leq i\leq N,
\end{equation}
where, for $1\leq j\leq N$, $(t_{jk}, k\geq 1)$  are $N$ i.i.d. Poisson point processes on $\R_+$ with rate $\lambda$ and, for $i$, $k\in\N$, the random variables $(W^j_{ik}, 1\leq j\leq N)$ are i.i.d. with the same distribution as $W_1$. 
\end{proposition}
\begin{proof}
The proof relies on a backward coupling argument, see Levin et al.~\cite{Levin} for a general presentation of so called coupling from the past methods and Loynes~\cite{Loynes} for one of its early uses.  Let $({\cal N}_j, 1\leq j\leq N)$ be $N$ i.i.d. Poisson Processes on $\R$ with rate $\lambda$, for $1\leq j\leq N$, ${\cal N}_j$ is the sequence of instants when the $j$th node spikes. Note that the time interval considered is $(-\infty, +\infty)$. Assume that for some fixed $T>0$, $X_i(-T)=0$ for all $1\leq j\leq N$, then, by using the invariance properties of Poisson processes, it is not difficult to see that $(X_i(0),1\leq i\leq N)$ has the same distribution of the state of the network at time $T$ when it starts empty at time $0$.

For $1\leq i\leq N$, if $t_{i,-1}$ is the last instant of spike of node $i$ before time $0$,  then if  $-T\leq t_{i,-1}$ state of this node at time $0$ is determined by the spikes of the other nodes  after  time $t_{i,-1}$.  If node $j\not= i$ spikes at time $s$, $t_{i,-1}\leq s\leq 0$, then the contribution at time $0$ for node $i$ is the value of the spike multiplied by $\exp(s)$.  Consequently, if $T$ is sufficiently large the value of $X(0)$ does not depend on $T$ and its distribution is the law of the vector given by Relation~\eqref{mi-const}. The proposition is proved. 
\end{proof}

\begin{proposition}
When the firing rate is constant and equal to $\lambda$, the Laplace transform of the state of a node at equilibrium   is given by, for $\xi\geq 0$, 
\[
\E\left(e^{-\xi X_1}\right)=\int_0^{+\infty} 
\exp\left(-\lambda(N-1)\int_0^x \left(1-\widetilde{W}\left(\xi e^{-u}\right)\right)\,\diff u\right)\lambda e^{-\lambda x}\,\diff x.
\]
where $\widetilde{W}(\xi)=\E(\exp(-\xi W))$ is the Laplace transform of $W$ at $\xi$.
\end{proposition}
This formula gives $\P(X_1{=}0){=}1/N$, which is simply the probability that the node is the last one which spiked. Similarly, the expected value at equilibrium is given by $\E(X_1)=(N-1)\E(W)\lambda/(\lambda+1)$
\begin{proof}
With the same notations as before, ${\cal M}={\cal N}_2+\cdots+{\cal N}_N=(s_n)$ is a Poisson process with rate $\lambda(N-1)$, and the above proposition gives that
\begin{align*}
\E\left(e^{-\xi X_1}\right)&=\E\left(\exp\left(-\xi \sum_{n\geq 1} W_{n}e^{-s_n}\ind{s_n\leq t_{11}}\right)\right)\\
&=\left.\E\left(\E\left(\exp\left(-\xi \sum_{n\geq 1} W_{n}e^{-s_n}\ind{s_n\leq t_{11}}\right)\right| (s_n), t_{11}\right)\right)
\end{align*}
where $(W_{1}(s), s\in{\cal M})$ are i.i.d. with the same distribution as $W$.  Consequently,  one obtains that
\begin{align*}
\E&\left(e^{-\xi X_1}\right)=\left.\E\left(\prod_{s_n\leq t_{11}}\E\left(\exp\left(-\xi  W_{1}e^{-s_n}\right)\right| s_n, t_{11}\right)\right)\\
&=\E\left(\prod_{s_n\leq t_{11}} \widetilde{W}\left(\xi e^{-s_n}\right)\right)
=\E\left( \exp\left(-\int_0^{t_{11}} g(u)\,{\cal M}(\diff u)\right)\right)
\end{align*}
with 
\[
g(u)\eqdef-\log\left(\widetilde{W}\left(\xi e^{-u}\right)\right).
\]
This gives the relation 
\[
\E\left(e^{-\xi X_1}\right)=\int_0^{+\infty} \E\left(\exp\left(-\int_0^x g(u)\,{\cal M}(\diff u)\right)\right)\lambda e^{-\lambda x}\,\diff x,
\]
since $t_{11}$ is exponentially distributed with parameter $\lambda$. 

The point process ${\cal M}$ being Poisson with rate $\lambda(N-1)$, from a classical formula for its Laplace transform, see Proposition~1.5 of Robert~\cite{Robert} for example, one gets 
\[
\E\left(\exp\left(-\int_0^x g(u)\,{\cal M}(\diff u)\right)\right)=
\exp\left(-\lambda(N-1)\int_0^x \left(1-e^{-g(u)}\right)\,\diff u\right).
\]
The Laplace transform of $X_1$ can thus be expressed as 
\[
\E\left(e^{-\xi X_1}\right)=\int_0^{+\infty} 
\exp\left(-\lambda(N-1)\int_0^x \left(1-\widetilde{W}\left(\xi e^{-u}\right)\right)\,\diff u\right)\lambda e^{-\lambda x}\,\diff x.
\]
The proposition is proved. 
\end{proof}
One concludes with a limiting regime which will be analyzed in a more general framework in the following. With little effort, it gives an idea of the results which can be obtained when the size of the network gets large, for example that the states of the nodes become independent in the limit. For this limiting regime, the size $N$ of the network goes to infinity and the rescaling is achieved through the  values of  spikes which are  of the order of $1/N$. 
The proposition shows in fact  the mean-field convergence of the invariant distribution of the state of the network.  See Sznitman~\cite{Sznitman} for an introduction on this topic.  
\begin{proposition}[A large network at equilibrium for constant firing rate]
When the firing rate is constant and equal to $\lambda$ and the value of a spike is $V_1/N$ for some integrable random variable $V_1$ and if  $(X_i^N)$ is the vector whose distribution is the equilibrium distribution of the state of the network then
\begin{enumerate}
\item the sequence of random variables $(X_1^N)$ converges in distribution to a random variable $X_1^\infty$ whose  distribution  has the density
\[
\frac{1}{\E(V_1)}\left(1-\frac{u}{\lambda\E(V_1)}\right)^{\lambda-1}\text{ for } u\in [0,\lambda \E(V_1)].  
\]
\item For fixed $i$ and $j$, $1\leq i<j\leq N$, the random variables $X_i^N$ and $X_j^N$ are asymptotically independent when $N$ gets large. 
\end{enumerate}
\end{proposition}
\begin{proof}
By symmetry, one can take $i=1$ and $j=2$, for $\ell=1$, $2$, if 
\[
Y_\ell^N\stackrel{\text{def}}{=}\sum_{j=3}^N \sum_{k\geq 1} \frac{V^j_{ik}}{N}e^{-t_{jk}} \ind{t_{jk}\leq t_{\ell 1}},\quad 1\leq i\leq N,
\]
then, as $N$ gets large,  the distribution of  $(X_1^N,X_2^N)$  is arbitrarily close to the distribution of $(Y_1^N,Y_2^N)$, just because the contribution of the spikes of node $1$ to the state of node $2$ are of the order of $1/N$ and vice versa. For $\xi_1$, $\xi_2\geq 0$, by using the same method as in the proof of the above proposition, one gets that
\begin{multline*}
\E\left(e^{-\xi_1Y_1^N-\xi_2Y_2^N}\right)=\\
\E\left(\exp\left[-\lambda(N{-}2)\int_0^{+\infty}\left(1{-}\widetilde{V}\left(\frac{\xi_1}{N}  e^{-u}\ind{u\leq t_{11}}\right)\widetilde{V}\left(\frac{\xi_2}{N} e^{-u}\ind{u\leq t_{21}}\right)\right)\,\diff u\right]\right),
\end{multline*}
where $\widetilde{V}$ denotes the Laplace transform of $V_1$. The equivalence $1-\widetilde{V}(x)\sim x\E(V_1)$ when $x$ goes to $0$ gives the relation
\[
\lim_{N\to+\infty} \E\left(e^{-\xi_1Y_1^N-\xi_2Y_2^N}\right)=H(\xi_1)H(\xi_2),
\]
with
\[
H(\xi)=\E\left(\exp\left[-\xi \lambda\E(V_1)\left(1- e^{-t_1}\right)\right]\right).
\]
This gives the asymptotic independence and the identification of the limit. The proposition is proved. 
\end{proof}

The analysis of the constant firing-rate model is instructive in many regards: it provides a completely solvable model for which no explosion of the firing rate is found. By comparison with the constant firing-rate case, coupling methods may allow to show that there is no explosion of the total firing-rate in the network (see~\cite[Appendix A.1]{demasi-presutti:14}), i.e. that the probability of occurrence of a large number of spikes in a fixed interval is small. We will come back to this property in the forthcoming section. Let us just state that this is an important from the biological viewpoint: consistently with the actual firing of neurons and in contrast with what happens for the stochastic integrate-and-fire neuron, there is no explosion of the firing rate and non-explosion of the membrane potential. 

%% file: McKeanVlasov.tex
\section{Analysis of the McKean Vlasov process} \label{McKeanSec}
We shall prove in particular in section~\ref{sec:Convergence} that  if the i.i.d. sequence $(W_{ij})$ has the same distribution as $(V_{ij}/N)$ where $(V_{ij})$ are i.i.d. with the same distribution as $V$, then 
$(X_1^N(t))$ converges in law towards the distribution of the stochastic process $(X(t))$ such that, for all $t\geq 0$, $b(X(t))$ is integrable and it satisfies the SDE:
\begin{multline}\label{n2}
\diff X(t)=\Big(\W \E(b(X(t))) -X(t)\Big)\diff t\\-X(t-)\int \ind{0\leq u\leq b(X(t-))} {\cal N}(\diff u,\diff t,\diff z),
\end{multline}

The object of this section is to show that this \emph{McKean-Vlasov} equation defines a unique process. Section~\ref{sec:Invariant} will characterize its stationary solutions. Throughout this section, we denote for $(U(t))$ a locally bounded process and $T>0$ the quantity $\|U\|_{T}$ defined as 
\[
\|U\|_{T}=\sup\{|U(t)|, 0\leq t\leq T\}.
\]
\begin{lemma}\label{lem-est}
If  $b$ is  a non-decreasing $C^1$-function on $\R_+$ and   ${\cal P}$ is a Poisson process on $\R_+^2$ with rate $1$, for any    non-negative  locally bounded Borelian function $(u(t))$   there exists a  unique solution  $(Z_u(x,t))$ of the  SDE
\begin{equation}\label{Vu}
\diff Z_u(x,t)=-Z_u(x,t)\diff t + u(t)\diff t -Z_u(x,t-) {\cal P}([0, b(Z_u(x,t-))],\diff t)
\end{equation}
with initial condition $x>0$.  For any couple of  non-negative  locally bounded Borelian functions $u$ and $v$ on $\R_+$, for $t\leq T$, the relation 
\begin{equation}\label{Est-Vu}
\E\left(\|Z_u-Z_v\|_{t}\right)\leq  e^{D_T t} \int_0^t \|u-v\|_{s}\,\diff s
\end{equation}
holds with  $D_T=1+(x+T\|u\|_{T})(1+\|b'\|_{ x+T\|u+v\|_{T}})$.
\end{lemma}
\begin{proof}
For a non-negative Borelian function $(u(t))$, the existence and uniqueness of a solution to the SDE~\eqref{Vu} is straightforward. 
Let  $u$ and $v$ be  non-negative locally bounded Borelian functions  on $\R_+$. Note that, almost surely, 
\[
\|Z_u\|_{T}\leq x+ \int_0^T u(s)\,\diff s\leq x+T\|u\|_{T}. 
\]
One has  to estimate, for $0\leq t\leq T$,
\[
\Delta(t)\eqdef\left|\int_0^t Z_{u}(s-) {\cal P}([0, b(Z_{u}(s-))],\diff s)-\int_0^t Z_{v}(s-) {\cal P}([0, b(Z_{v}(s-))],\diff s)\right|
\]
then
\begin{multline*}
\Delta(t)\leq \int_0^t \left|Z_{u}(s-)-Z_{v}(s-)\right| {\cal P}([0, b(Z_{u}(s-))],\diff s)\\+\left|\int_0^t Z_{v}(u-) {\cal P}([b(Z_{u}(s-)), b(Z_{v}(s-))],\diff s)\right|
\end{multline*}
hence
\begin{align*}
\E(\|\Delta\|_{T})&\leq \E\left(\int_0^t \|Z_{u}-Z_{v}\|_{s} {\cal P}([0, b(Z_{u}(s-))]\,\diff s)\right)\\&\hspace{2cm}+\E\left(\left|\int_0^t Z_{v}(u-) {\cal P}([b(Z_{u}(s-))\, b(Z_{v}(s-))]\,\diff s)\right|\right)\\
&\leq b(x+T\|u\|_{T})\int_0^t \E(\|Z_{u}-Z_{v}\|_{s}) \diff s\\
&\hspace{2cm}+(x+T\|u\|_{T})\int_0^t \E(|b(Z_{v}(s))-b(Z_{u}(s))|)\,\diff s.
\end{align*}
The SDE associated to $u$ and $v$ give, for $0\leq t\leq T$ and a convenient constant $D_T$,
\[
\E(\|Z_{u}-Z_{v}\|_{t})\leq D_T \int_0^t\E(\|Z_{u}-Z_{v}\|_{s})\,\diff s+ \int_0^t\|u-v\|_{s}\,\diff s,
\]
with $D_T$ as defined above.
Gronwall's Lemma completes the proof of the lemma. 
\end{proof}

\begin{lemma}[Integrability]\label{lem:integrability}
If $b$ is  a non-decreasing unbounded $C^1$ function such that there exist  $\gamma>0$ with $3\gamma \W<1$ and $c > 0$  and, for all $x\geq 0$<  $$b'(x)\leq \gamma b(x) + c$$  and, for any $0\leq C\leq +\infty$,   if $(Z^C(t))$ is a solution of the following SDE
\begin{multline}\label{eqZC}
\diff Z^C(t)=\Big(\W [C \wedge\E(b(Z^C(t)))] -Z^C(t)\Big)\diff t\\-Z^C(t-)\int \ind{0\leq u\leq b(Z^C(t-))} {\cal N}(\diff u,\diff t,\diff z),
\end{multline}
 with an  initial condition $Z(0)$ independent of $C$ such that $Z(0)$ and $b^3(Z(0))$ are integrable 
then, for $p\in \{1,2,3\}$, 
\[
\sup_{\substack{t\geq 0\\0\leq C\leq +\infty}} \E(Z^C(t)) <+\infty\quad \text{ and }\quad \sup_{\substack{t\geq 0\\0\leq C\leq +\infty}} \E\left(b\left(Z^C(t)\right)^p\right) <+\infty.
\]
\end{lemma}
\begin{proof}
 For a fixed $t$,  it is easily seen that the non-negative random variable $Z^C(t)$ is integrable.  Define $\mu(t)=\E(Z^C(t))$. From Equation~\eqref{eqZC}, one gets
\begin{multline*}
\mu(t)\leq \mu(0) + \int_0^t\left[\rule{0mm}{4mm} -\mu(s) + \W\E(b(Z^C(s))) - \E(Z^C(s) b(Z^C(s)))\right]\,\diff s\\ = \mu(0) + \int_0^t -\mu(s) + \E(\Phi(Z^C(s)))\,\diff s.
\end{multline*}
with $\Phi(x)= (\W -x) b(x)$. The equivalence $\Phi(x)\sim -xb(x)$ as $x$ gets large and Gronwall's lemma give the   boundedness of the first moment. 
The proof of lemma~\ref{lemelem} in the Appendix shows that, there exists some $\eps>0$ such that, for $1\leq p\leq 3+\eps$,  the derivative of $b^p$ is also upper-bounded by $\gamma_1 b^p(x)+c$ with $\gamma_1<1/\W$. It is thus enough to prove the boundedness of the first moment $B(t)=\E(b(Z^C(t))$ assuming $\gamma<1/\W$. Equation~\eqref{eqZC} gives
	\begin{multline*}
B(t)\leq B_0 + \int_0^t\left[\rule{0mm}{4mm}  \E(b'(Z^C(s))(-Z^C(s)+\W B(s))\right.\\\left. \rule{0mm}{4mm}  + (b(0)-b(Z^C(s)))b(Z^C(s)))\right] \,\diff s.
\end{multline*}
	and by Cauchy-Schwarz' inequality using the fact that $\gamma \W<1$,
	\begin{align*}
		B(t) &\leq B(0) + \int_0^t \left(\E\left[\rule{0mm}{4mm} {-}b'(Z^C(s))Z^C(s){+}(b(0)-b(Z^C(s)))b(Z^C(s))\right] \right.\\ & \left. \hspace{6cm}\hfill{+} \E(b'(Z^C(s)))\W B(s) \rule{0mm}{4mm} \right) \,\diff s\\
		&\leq B(0) + \int_0^t \E((b(0)-b(Z^C(s)) + \gamma \W)b(Z^C(s)) + \gamma \W b(Z^C(s))^2) \,\diff s\\
		& \leq B(0) + \int_0^t (b(0)+\W \gamma) B(s)  + (\gamma \W-1) B(s)^2) \,\diff s.
	\end{align*}
	We conclude using lemma~\ref{elmGrow}, clearly the associated upper bound for $B(t))$ does not depend on $C$. 
\end{proof}

We can now state our main result on the existence and uniqueness of a solution to mean-field  equations~\eqref{n2}. 
\begin{theorem}[Existence and Uniqueness of the McKean-Vlasov Process]\label{McVTheo}
If $b$ is  a non-decreasing unbounded $C^1$ function and if there exist  $\gamma\geq 0$ with $3\gamma \W<1$ and $c > 0$ such that, for all $x\geq 0$, 
\[
b'(x)\leq \gamma b(x) + c,
\]
 then for any $T>0$, there exists a unique c\`adl\`ag process $(Z(t))$ satisfying the stochastic differential equation
\begin{equation}\label{Vlasov}
\diff Z(t)=-Z(t)\diff t +\W \E(b(Z(t)))\diff t
-Z(t-) {\cal P}([0, b(Z(t-))],\diff t),
\end{equation}
and with initial condition $x>0$. 
\end{theorem}
\begin{proof}
For $C>0$ we define by induction the sequence of processes $(Z_n(t))$ by $Z_0(t)= x\exp(-t)$ and for $t>0$ and   $n\geq 1$,
\[
\diff Z_n(t)=-Z_{n}(t)\diff t +\W  \big[\E(b(Z_{n-1}(t)))\wedge C \big]\diff t -Z_{n}(t-) {\cal P}([0, b(Z_{n}(t-))],\diff t),
\]
with $Z_n(0)=x$. It is easy to show that we have:
\[
Z_{n}(t)\leq  xe^{-t}+C (1-e^{-t})\leq x+C.
\]
For $T>0$, Lemma~\ref{lem-est} and the above relation show that there exists a constant $D_T$ independent of $n$ such for $0\leq t\leq T$,
\[
\E\left(\|Z_{n+1}-Z_{n}\|_{t}\right)\leq D_T\int_0^t \|u_{n}-u_{n-1}\|_{s}\,\diff s,
\]
with $u_n(t)=\E(b(Z_n(t)))$. This implies that:
\[
\E\left(\|Z_{n+1}-Z_{n}\|_{t}\right)\leq D_T\int_0^t \E\left(\|b(Z_{n})-b(Z_{n-1})\|_{s}\right)\,\diff s,
\]
and thanks to the deterministic bound on the sequence of processes $(Z_n)$, we have:
\[
\E\left(\|Z_{n+1}-Z_{n}\|_{t}\right)\leq K_T\int_0^t \E\left(\|Z_{n}-Z_{n-1}\|_{s}\right)\,\diff s,
\]
with $K_T=D_T \|b'\|_{x+C}$.

We hence have:
\[
\E\left(\|Z_{n+1}-Z_{n}\|_{t}\right)\leq \frac{(K_Tt)^n}{n!} \|z(x,\cdot)\|_{T}
\]
From this relation  one gets that 1) the sequence of processes $(Z_N(t), 0\leq t\leq T)$ is converging almost surely uniformly on compact sets to a solution $(Z(t), 0\leq t\leq T)$ of SDE~\eqref{Vlasov}.
If $(\widetilde{Z}(t), 0\leq t\leq T)$ is another solution of this SDE starting from $x$, then necessarily  $\widetilde{Z}(t)\leq z(t)$ for all $0\leq t <T$, consequently the relation
\[
\E\left(\|Z_{n+1}-\widetilde{Z}\|_{t}\right)\leq K_T\int_0^t \E\left(\|Z_{n}-\widetilde{Z}\|_{s}\right)\,\diff s,
\]
holds, hence
\[
\E\left(\|Z-\widetilde{Z}\|_{t}\right)\leq K_T\int_0^t \E\left(\|Z-\widetilde{Z}\|_{s}\right)\,\diff s,
\]
so $Z$ and $\widetilde{Z}$ are identical. We have therefore shown that there exists a unique solution to the equation:
\begin{multline}\label{eq:Cutoff}
		\diff Z^C(t)=-Z^C(t)\diff t +\W  \big[\E(b(Z^C(t)))\wedge C \big]\diff t\\ -Z^C(t-) {\cal P}([0, b(Z^C(t-))],\diff t)
\end{multline}
with $Z(0)=x$. 
Lemma~\ref{lem:integrability} actually ensures that there exists some constant $C_0$ such that  $\E(b(Z^C(t)))<C_0$ for all $t\geq 0$ and $C>0$. Hence, for $C\geq C_0$, $b(Z^C(t))$ is integrable and  $(Z^C(t))$ is a solution of~\eqref{Vlasov}. Conversely, by using again Lemma~\ref{lem:integrability} with $C=+\infty$, any solution of Equation~\eqref{n2}  is also solution of Equation~\eqref{eq:Cutoff}. The theorem is proved. 
\end{proof}

%% file: Mean-Field.tex
\section{Mean-Field Asymptotics}\label{sec:Convergence}
In this section, one considers the asymptotic regime of these networks when  the number $N$ of nodes of the network goes to infinity. The interaction between nodes is as follows: for $1\leq i\not=j\leq N$,  when node $i$ fires the value of the state of node $j$ is increased by $W_{ij}^N=V_{ij}/N$.  The variables $(V_{ij})$ are i.i.d.  integrable random variables with distribution $V(\diff x))$, with a slight abuse of notation $V$ denotes in the following a random variable with such a distribution.
The associated Markov process is denoted by $(X_i^N(t))$. We recall the SDE equations~\eqref{SDE1} in this context, for $1\leq i\leq N$,  one has
\begin{multline}\label{SDEMF}
\diff X^N_i(t)=-X^N_i(t)\diff t +\frac{1}{N}\sum_{j\not=i} \int_{\R_+^2} z_i\ind{0\leq u\leq b(X^N_j(t-))} {\cal N}_j(\diff u,\diff z,\diff t)
\\-X_i^N(t-)\int_{\R_+^2} \ind{0\leq u\leq b(X^N_i(t-))} {\cal N}_i(\diff u,\diff z,\diff t),
\end{multline}
where ${\cal N}_i$ is a Poisson processes with intensity measure  given by $\diff u\otimes V(\diff z)\otimes \diff t$. The Poisson processes $({\cal N}_j, 1\leq j\leq N)$ are independent 

The  {\em empirical distribution}  is denoted by  $(\Lambda_N(t))$, for any continuous function $\phi$ on $\R_+$, 
\[
\croc{\Lambda_N(t),\phi}=\frac{1}{N}\sum_{i=1}^N \phi(X_i^N(t)). 
\]
The main result of this section is that  under appropriate conditions a mean-field convergence holds:  The sequence $(\Lambda_N(t))$ of random measures valued processes converges in distribution to the distribution of the McKean-Vlasov process analyzed in Section~\ref{McKeanSec}. 

The strategy of the proof is the following:  first it is shown  that the scaled moment of the total firing rate of the network
\[
\frac{1}{N}\sum_{i=1}^N b(X^N_i(t))
\]
is, with high probability, bounded uniformly on any finite time interval.   Then,  by using the stochastic  evolution equations of $(\Lambda_N(t))$, it is then proved  that  for any  continuous function $\phi$ on $\R_+$ with compact support,  the sequence of processes $(\croc{\Lambda_N(t),\phi})$ is tight for the topology of the uniform norm, in particular any of its limiting points is a continuous process. One concludes by a uniqueness result proved in the appendix. 

\medskip
Concerning the main parameters and the initial state  of the network, the assumptions are given below.

\renewcommand{\labelenumi}{(\alph{enumi})}
{\bf Assumptions MF}
\begin{enumerate}
\item Growth Condition.\\
{\em The firing rate function $x\mapsto b(x)$ is assumed to be  $C^1$, non-decreasing and such that there exist  $3\E(V)\gamma < 1$ and $c>0$  such that
  \begin{equation}\label{eq:assumption}
    b'(x) \leq \gamma b(x)+c
  \end{equation}
holds for any $x\geq 0$.}
\item Bounded Support.\\
{\em The distribution $V(\diff x))$ has a bounded support,  there exists some $S_V>0$ such that $V([0,S_V])=1$. }
\item {Initial Conditions.}\\
{\em The random variables $(X_i^N(0), 1\leq i\leq N)$ are i.i.d. with law $m_0$ having a bounded support.}
\end{enumerate}
\renewcommand{\labelenumi}{(\arabic{enumi})}
Assumption~(MF-a) implies that, for any $a>0$, the ratio
${b(x+a)}/{b(x)}$ 
is bounded as $x$ gets large, i.e. a slow growth at infinity. Note that polynomial functions satisfy this assumption. See the proof of Lemma~\ref{lemelem} in the Appendix.  Additionally, for convenience, it will be assumed that $b(0)=0$ in the following. It turns out that the case $b(0)>0$ is easier from the point of view of the mean-field analysis of this section. Indeed, in this case, the nodes are ``refreshed'' at a minimal positive rate, in this way there is a maximal, state independent, interval between two spikes of a given node. 

\subsection{Stochastic Evolution Equations for the Empirical Distribution}
Let $f$ a $C^1$-function on $\R_+$  then, from Equation~\eqref{SDEMF} and Proposition~\ref{PoisMart} of the Appendix, one gets that, for $1\leq i\leq N$, 
\begin{multline}\label{SDE2}
f(X_i^N(t))=f(X_i^N(0))-\int_0^t X_i^N(u)f'(X_i^N(u))\,\diff u\\ + \sum_{j\not=i} \int_0^t \int_{\R_+}\left(f\left(X_i^N(u)+\frac{v}{N}\right)-f(X_i^N(u))\right) b(X_j^N(u))\,\diff u V(\diff v)
\\ +\int_0^t\left[f(0)-f(X^N_i(u))\right]b(X^N_i(u))\,\diff u+M^N_{f,i}(t),
\end{multline}
where $(M^N_{f,i}(t))$ is the local martingale defined by
\begin{multline*}
\int_{s{=}0}^t\! \int_{\R_+^2}\!\! \left[\rule{0mm}{4mm}f(0){-}f\!\left(X^N_i(s-)\right)\right] \!\!\left[\ind{0\leq u\leq b(X^N_j(s-))}  {\cal N}_j(\diff u,\diff z,\diff s){-} b(X^N_j(s))\,\diff s  \right]\\
+\sum_{j\not=i}\int_{s=0}^t  \int_{\R_+^2} \left(f\left(X_i^N(s-)+\frac{z_i}{N}\right)-f(X_i^N(s-))\right)\\
\left[
\ind{0\leq u\leq b(X^N_j(s-))} {\cal N}_j(\diff u,\diff z,\diff s)
-b(X^N_j(s))\,\diff s V(\diff z_i)
\right].
\end{multline*}

Provided that the local martingales  $(M^N_{f,i})$, $i=1$, \ldots,$N$ are locally square integrable, the associated previsible increasing processes are given by
\begin{multline}\label{crocM}
\croc{M^N_{f,i}}(t)= \int_{0}^t \left[f(0)- f(X^N_i(s)\right]^2b(X^N_i(s))\,\diff s \\
+\sum_{j\not=i}\int_0^t  b(X^N_j(s))\,\diff s  \int_{\R_+} \left[f\left(X_i^N(s)+\frac{v}{N}\right)-f\left(X_i^N(s)\right)\right]^2\,V(\diff v),
\end{multline}
and, for $1\leq i<j\leq N$,
\begin{multline}\label{crocMM}
\croc{M^N_{f,i},M^N_{f,j}}(t)=\\ \int_{0}^t \left[f(0)- f(X^N_i(s)\right]\left[f\left(X_j^N(s)+\frac{v}{N}\right)-f\left(X_j^N(s)\right)\right]b\left(X^N_i(s)\right)\,\diff s \\
\int_{0}^t \left[f(0)- f(X^N_j(s)\right]\left[f\left(X_i^N(s)+\frac{v}{N}\right)-f\left(X_i^N(s)\right)\right]b\left(X^N_j(s)\right)\,\diff s \\
+\sum_{k\not\in\{i,j\}}\int_0^t   \int_{\R_+} \left[f\left(X_i^N(s)+\frac{v}{N}\right)-f\left(X_i^N(s)\right)\right]\,V(\diff v)
\\ \times  \int_{\R_+} \left[f\left(X_j^N(s)+\frac{v}{N}\right)-f\left(X_j^N(s)\right)\right]\,V(\diff v) \, b\left(X^N_k(s)\right)\diff s
\end{multline}
by Proposition~\ref{PoisMart} of the Appendix. Equation~\eqref{SDE2} gives therefore  the following relation for the empirical measure
\begin{multline}\label{EmpSDE}
\croc{\Lambda_N(t),f}=\croc{\Lambda_N(0),f}-\int_0^t \int_{\R_+}\croc{\Lambda_N(u), \odot f'(\cdot)}\,\diff u\\
+N\int_0^t \int_{\R_+}\croc{\Lambda_N(u),\left(f\left(\cdot+\frac{v}{N}\right)-f(\cdot)\right)} \croc{\Lambda_N(u),b}\,\diff u V(\diff v)\\
-\int_0^t\croc{\Lambda_N(u),\left(f\left(\cdot+\frac{v}{N}\right)-f(\cdot)\right)b(\cdot)}\,\diff u V(\diff v)\\-\int_0^t\croc{\Lambda_N(u),(f(0)-f(\cdot))b(\cdot)}\,\diff u+M_{f}^N(t),
\end{multline}
where  $(M_{f}^N(t))$ is the martingale
\begin{equation}\label{martMf}
M_{f}^N(t)=\frac{1}{N}\sum_{i=1}^N M_{f,i}^N(t),
\end{equation}
The corresponding  previsible increasing process is given by
\begin{equation}\label{crocMf}
\croc{M_{f}^N}(t)=\frac{1}{N^2}\left(\sum_{i=1}^N \croc{M^N_{f,i}}(t)+2\sum_{i<j} \croc{M^N_{f,i},M^N_{f,j}}(t)\right).
\end{equation}
\subsection{Estimates for the scaled firing rate}
In this section it is proved that the scaled  moment of the firing rate 
\[
\croc{\Lambda_N(t),b}=\frac{1}{N}\sum_{i=1}^N b\left(X_i^N(t)\right),
\] 
remains with high probability within a finite interval. One starts with a result on the boundedness of some of its moments.
\begin{lemma}\label{EstLem}
Under Assumptions~(MF), there exists $\delta>3$ such that relation
\[
\sup_{N\geq 1} \sup_{t\geq 0} \E(\croc{\Lambda_N(t),b^\delta}) <+\infty, 
\]
holds.
\end{lemma}
\begin{proof}
By Lemma~\ref{lemelem} of the Appendix and Assumptions~(MF) there exists $\delta>3$ such that 
\begin{equation}\label{eq21}
b^\delta(x+a)-b^\delta(x) \leq a \left(\gamma_1 b^\delta(x)+c_1 \right), \forall a\in(0,\eta_b) \text{ and } \forall x\geq 0,
\end{equation}
for some $\eta_b>0$ and with $\gamma_1<1/\E(V)$. 

For $K>0$. let $\tau_K=\inf\{t\geq 0: \croc{\Lambda_N(t),b^{\delta+1}} \geq K\}$. 
Holder's Inequality shows that  for all $t\geq 0$ the random variables $\croc{\Lambda_N(t\wedge\tau_K),b^p}$, $1\leq p\leq \delta+1$, are integrable. 
By  taking $f=b^\delta$ in Equation~\eqref{EmpSDE}, the optional stopping theorem gives  the relation
\begin{multline*}
\E\left(\croc{\Lambda_N(t\wedge\tau_K),b^\delta}\right)\leq \croc{\Lambda_N(0),b^\delta}
\\+N\int_0^t \int_{\R_+}\E\left(\croc{\Lambda_N(u\wedge\tau_K),b^\delta\left(\cdot+\frac{v}{N}\right)-b^\delta(\cdot)} \croc{\Lambda_N(u\wedge\tau_K),b}\right)\,\diff u V(\diff v)
\\ -\int_0^t\E(\left(\croc{\Lambda_N(u\wedge\tau_K),b^{1+\delta}})\right)\,\diff u,
\end{multline*}
Recall that $S_V$ is an upper bound for the support  of the distribution $V$,  $N$ is chosen sufficiently large so that  $S_V/N\leq \eta_b$, 
from  Inequality~\eqref{eq21},  one gets
\begin{multline*}
\E\left(\croc{\Lambda_N(t\wedge\tau_K),b^\delta}\right)\leq  \croc{\Lambda_N(0),b^\delta}\\
+\gamma_1\E(V)\int_0^t\E\left(\croc{\Lambda_N(u\wedge\tau_K),b^\delta}\croc{\Lambda_N(u\wedge\tau_K),b}\right)\,\diff u \\
+ c_1\E(V)\int_0^t \E\left(\croc{\Lambda_N(u\wedge\tau_K),b}\right)\,\diff u -\int_0^t\E(\left(\croc{\Lambda_N(u\wedge\tau_K),b^{1+\delta}})\right)\,\diff u. 
\end{multline*}
Holder's Inequality and the fact that $\Lambda_N$ is a probability distribution give
\begin{align*}
\croc{\Lambda_N(u{\wedge}\tau_K),b}&\leq \croc{\Lambda_N(u{\wedge}\tau_K),b^\delta}^{1/\delta} \\
\croc{\Lambda_N(u{\wedge}\tau_K),b^{1+\delta}}&\geq\croc{\Lambda_N(u{\wedge}\tau_K),b^\delta}^{(1+\delta)/\delta}
\end{align*}
 hence
\begin{multline*}
\E\left(\croc{\Lambda_N(t\wedge\tau_K),b^\delta}\right)\leq  C_0
-(1-\gamma_1\E(V))\int_0^t\E\left(\croc{\Lambda_N(u\wedge\tau_K),b^\delta}^{(1+\delta)/\delta}\right)\,\diff u \\
+c_1\E(V)\int_0^t \E\left(\croc{\Lambda_N(u\wedge\tau_K),b^\delta}^{1/\delta}\right)\,\diff u,
\end{multline*}
where 
\[
C_0\stackrel{\text{def.}}{=}\sup_{N\geq 1}\E(\croc{\Lambda_N(0),b^\delta}).
\]
With  again Holder's Inequality and the fact that $\gamma_1\E(V)<1$, one gets finally
 \begin{multline*}
\E\left(\croc{\Lambda_N(t\wedge\tau_K),b^\delta}\right)\leq  C_0
-(1-\gamma_1\E(V))\int_0^t\left[\E\left(\croc{\Lambda_N(u\wedge\tau_K),b^\delta}\right)\right]^{\delta/(1+\delta)}\,\diff u \\
+c_1\E(V)\int_0^t \left[\E\left(\croc{\Lambda_N(u\wedge\tau_K),b^\delta}\right)\right]^{1/\delta}\,\diff u,
\end{multline*}
By using the inequality $\gamma_1\E(V)<1$ and, since $\delta/(1+\delta)>1/\delta$,  Proposition~\ref{elmGrow} of the Appendix, one gets that there exists a finite constant $C_0$ independent of $K$ such that
\[
\sup_{N\geq 1}\sup_{t\geq 0}\E\left(\croc{\Lambda_N(t\wedge\tau_K),b^\delta}\right)\leq C_0,
\]
on concludes the proof by letting $K$ go to infinity. 
\end{proof}
\begin{proposition}[Control of the scaled firing rate]\label{EstProp}
Under Assumptions~(MF), for any $T>0$ there exists $\kappa>1$ and some constant $C_0$ such that the 
\begin{equation}\label{estD}
\lim_{N\to+\infty}\P\left(\sup_{0\leq t\leq T}\croc{\Lambda_N(t),b^\kappa} \geq C_0\right)=0. 
\end{equation}
\end{proposition}
\begin{proof}
For simplicity, the proof is done for $\kappa=1$. The case $\kappa>1$ follows the same lines together with the same method as in the proof of the previous lemma with $\delta>3$. One first shows that there exists some constant $C_T$ such that, 
\begin{equation}\label{eqa0}
\sup_{0\leq t\leq T}\E\left(M_{b}^N(t)^2\right)=\E\left(\croc{M_{b}^N}(t)\right)\leq \frac{C_T}{N}. 
\end{equation}
The relations~\eqref{crocM}, \eqref{crocMM} and~\eqref{crocM}  are used in the case $f=b$. 
The first term of the right hand side of Relation~\eqref{crocM} for $\E(\langle M_{b}^N\rangle(t))$ gives the contribution
 \begin{equation}\label{eqa1}
\frac{1}{N^2}\sum_{i=1}^N\int_{0}^t\E\left(b(X^N_i(u))^3\right) \,\diff u=\frac{1}{N}\int_0^t\E\left(\croc{\Lambda_N(u),b^3}\right)\,\diff u
\end{equation}
and, for the second term,
\[
\frac{1}{N^2}\sum_{ 1 \leq  i\not=j\leq N}\int_0^t \int_{\R_+} \E\left[ b(X^N_j(u))  \left(b\left(X_i^N(u)+\frac{v}{N}\right)-b(X_i^N(u))\right)^2\right]\,\diff u V(\diff v).
\]
As before, if  $N$ is sufficiently large so that  $S_V/N\leq \eta_b$, Equation~\eqref{eq:Discrete} of the Appendix gives that this last term is upper bounded by
\begin{equation}\label{eqa2}
\frac{\E(V^2)}{N}\int_0^t \gamma_1\E\left(\croc{\Lambda_N(u),b^3}\right)+2\gamma_1c_1\E\left(\croc{\Lambda_N(u),b^2}\right) +c_1\E\left(\croc{\Lambda_N(u),b}\right)\,\diff u .
\end{equation}
With similar arguments, analogous bounds can be obtained for the second term  of Relation~\eqref{crocMM} for $\E(\langle M_{b}^N\rangle(t))$ involving the $\E(\langle M_{b,i}^N, M_{b,i}^N\rangle(t))$. Holder's Inequality and Lemma~\ref{EstLem} show that there exists a constant $C_T$ such that the upper bound~\eqref{eqa0} holds. Define
\[
M_{b}^{N,*}(T)\stackrel{\text{def.}}{=}\sup_{0\leq s\leq T}M_{b}^{N}(s),
\]
Doob's Inequality shows therefore that, for any $\eps>0$, there exists $N_0$ such that if $N\geq N_0$ then $\P(M_{b^2}^{N,*}(T)>1)\leq \eps$.
With  $f=b$ in Equation~\eqref{EmpSDE}, one gets that, for $0\leq t\leq T$,
\begin{multline}\label{eqa3}
\croc{\Lambda_N(t),b}\leq \croc{\Lambda_N(0),b}+M_{b}^{N,*}(T)
\\+N\int_0^t \int_{\R_+}\croc{\Lambda_N(u),b\left(\cdot+\frac{v}{N}\right)-b(\cdot)} \croc{\Lambda_N(u),b}\,\diff u V(\diff v)
\\ -\int_0^t\croc{\Lambda_N(u),b^2})\,\diff u.
\end{multline}
By the integrability condition of $b$ with respect to $m_0$ in Assumptions~(MF),  if $N$ is sufficiently large then, by the law of large numbers, 
\[
\P\left[\croc{\Lambda_N(0),b} > 1+\E\left[b(X_1^1)\right]\right]<\eps. 
\]
Let $ C_0\stackrel{\text{def.}}{=}2+\E(b(X_1^1(0)))\}$,
On the event $\{M_{b^2}^{N,*}(t)\leq 1\}$,  by using again Equation~\eqref{eq:Discrete} of the Appendix and Relation~\eqref{eqa0}, Relation~\eqref{eqa3} gives the inequality
\begin{multline}\label{eqa4}
\croc{\Lambda_N(t),b}\leq C_0
+(\gamma_1\E(V)-1)\int_0^t \int_{\R_+} \croc{\Lambda_N(u),b}^{2}\,\diff u V(\diff v)
\\+c_1\int_0^t\croc{\Lambda_N(u),b}\,\diff u.
\end{multline}
Proposition~\ref{elmGrow} allows to conclude that, on the event $\{M_{b^2}^{N,*}(t)\leq 1\}$, the random variable $\sup_{0\leq t\leq T}\croc{\Lambda_N(t),b}$ is uniformly bounded for all $N\geq 0$ (i.e., by a constant independent of $N$), proving~\eqref{estD}.
\end{proof}
\subsection{Tightness of $(\Lambda_N(t))$}
To prove this result, Theorem~3.7.1 of Dawson~\cite{Dawson:16} shows that it is enough to prove that, for any  continuous function $\phi$ on $\R_+$ with compact support,  the sequence of processes $(\croc{\Lambda_N(t),\phi})$ is tight for the topology of the uniform norm on compact sets. 
\begin{proposition}
The sequence of measure-valued processes $(\croc{\Lambda_N(t)})$ is tight for the convergence in distribution with continuous limits. 
\end{proposition}
\begin{proof}
Take $\phi$ a $C_1$-function  on $\R_+$ with compact support. By using the same method as in the proof of Proposition~\ref{EstProp}, one has that $\E\left(M_{\phi}^N(t)^2\right)$ is a $O(1/N)$, in particular $M_{\phi}^N(t)$ converges to $0$ in distribution for the uniform  norm on finite time interval. 

The modulus of continuity of $\croc{\Lambda_N,\phi}$ is defined as
\[
w_{\Lambda_N,\phi}(\delta)= \sup_{\substack{0\leq s\leq s'\leq t\\|s-s'|\leq \delta}} \left| \croc{\Lambda_N(s),\phi}-\croc{\Lambda_N(s'),\phi}\right|.
\]
For $s\leq s'$, 
\begin{multline*}
\croc{\Lambda_N(s),\phi}-\croc{\Lambda_N(s'),\phi}=-\int_s^{s'} \croc{\Lambda_N(u), \cdot\,\phi'(\cdot)}\,\diff u\\
+N\int_s^{s'} \int_{\R_+}\croc{\Lambda_N(u),\left(\phi\left(\cdot+\frac{v}{N}\right)-\phi(\cdot)\right)} \croc{\Lambda_N(u),b}\,\diff u V(\diff v)\\
-\int_s^{s'}\left(\croc{\Lambda_N(u),\phi}\right)\,\diff u+M_{\phi}^N(s)-M_{\phi}^N(s'),
\end{multline*}
take $\eta>0$ and $\eps>0$, by Equation~\eqref{estD}, there exists some $C_1$ such that, for all $N\geq 1$, 
\begin{equation}\label{aux2}
\P\left(w_{\Lambda_N,\phi}(\delta)>\eta\right)\leq \eps+\P\left(w_{\Lambda_N,\phi}(\delta)>\eta, \sup_{t\geq 0}\croc{\Lambda_N(t),b} \leq C_1\right).
\end{equation}
Now,  choose $\delta$ such that $\delta<\eta\min(1,(\E(W)C_1)^{-1})/(4K)$, then 
\[
\P\left(w_{\Lambda_N,\phi}(\delta)>\eta, \sup_{t\geq 0}\croc{\Lambda_N(t),b} \leq C_1\right)
\leq 
\P\left(\sup_{\substack{0\leq s\leq s'\leq t\\|s-s'|\leq \delta}} \left|M_{\phi}^N(s)-M_{\phi}^N(s') \right|\geq \eta/4\right),
\]
due to the convergence in distribution to $0$ of the martingale, this term can be made arbitrarily small for $N$ large. 

Hence, one has shown that the sequence of processes $(\croc{\Lambda_N(t),\phi})$ is tight for convergence in distribution for the uniform norm on compact sets. 
\end{proof}

Theorem~3.7.1 of Dawson~\cite{Dawson:16} shows that the sequence of measure-valued processes $(\Lambda_N(t))$ is tight, let $\Lambda(t)$ be one the limit of a given convergent subsequence $(\Lambda_{N_k}(t))$. In particular, for any continuous function with compact support, one has the convergence of the processes
\[
\lim_{N\to+\infty} \left(\croc{\Lambda_{N_k}(t),\phi}\right) = \left(\croc{\Lambda(t),\phi}\right).
\]
The next result shows that this convergence also for the function $(b(x))$. 
\begin{lemma}\label{lemb}
If $((\Lambda_{N_k}(t)))$ is a converging subsequence with $(\Lambda(t))$ as a limit, then 
\[
\lim_{k\to+\infty} \left(\croc{\Lambda_{N_k}(t),b}\right)=\croc{\Lambda(t),b}. 
\]
for the convergence in distribution. 
\end{lemma}
\begin{proof}
We first prove that the sequence of processes $(\croc{\Lambda_N(t),b})$ is tight. By a similar argument as before, for $0\leq s\leq s'\leq t$,
\begin{multline*}
\croc{\Lambda_N(s),b}-\croc{\Lambda_N(s'),b}=-\int_s^{s'} \croc{\Lambda_N(u), \cdot\,b'(\cdot)}\,\diff u\\
+N\int_s^{s'} \int_{\R_+}\croc{\Lambda_N(u),\left(b\left(\cdot+\frac{v}{N}\right)-b(\cdot)\right)} \croc{\Lambda_N(u),b}\,\diff u\,V(\diff v)\\
-\int_s^{s'}\left(\croc{\Lambda_N(u),b}\right)\,\diff u+M_{\phi}^N(s)-M_{\phi}^N(s'). 
\end{multline*}
 From the estimate~\eqref{estD} of Proposition~\ref{EstProp} one concludes that one can choose $\delta$ so that the left hand side of the above relation is arbitrarily small. The tightness has been proved. One has to identify the limit.

For $K>0$ there exists a function $\phi_K\in C_c(\R_+)$ which coincides with $b$  on the set $S_K=\{x:b(x)\leq K\}$ and such that $\phi_K\leq b$ and $K\mapsto \phi_K$ is increasing. We know that, for $t>0$ one has, for the convergence in distribution,
\[
\lim_{k\to+\infty}\croc{\Lambda_{N_k}(t),\phi_K}=\croc{\Lambda(t),\phi_K},
\]
and 
\begin{multline}\label{aux3}
\left|\croc{\Lambda_{N_k}(t),b}-\croc{\Lambda_{N_k}(t),\phi_K}\right|
\leq \croc{\Lambda_{N_k}(t),|b-\phi_K|\ind{b\geq K}}\\
\leq \croc{\Lambda_{N_k}(t),b\ind{b\geq K}}
\leq \frac{1}{K^{\kappa-1}}\croc{\Lambda_{N_k}(t),b^{\kappa}}.
\end{multline}
Proposition~\ref{EstProp} shows that the sequences of variables $(\croc{\Lambda_{N_k}(t),b^\kappa})$ is uniformly bounded on finite time intervals with high probability. Consequently, for $K$ sufficiently large, the left hand side of Relation~\ref{aux3} is arbitrarily close to $0$ in distribution for all $k\geq 1$. 

For $t\geq 0$, the monotone convergence theorem gives that
\[
\lim_{K\to+\infty} \croc{\Lambda(t),\phi_K}=\croc{\Lambda(t),b}.
\]
The convergence in distribution of $\croc{\Lambda_{N_k}(t),b}$ proved for a fixed $t$ is clearly valid for a vector of finite time marginals. By tightness one concludes that the process $(\croc{\Lambda_{N_k}(t),b})$ converges in distribution to $(\croc{\Lambda(t),b})$ for the topology of the uniform convergence on compact sets.
\end{proof}
\begin{theorem}[Mean-Field Convergence]\label{ThMF}
Under Conditions~(MF)   the sequence of processes $(X_1^N(t))$ converges in  distribution to the law of the unique process  $(Z(t))$  with initial distribution $m_0$ and solution of  the SDE~\eqref{Vlasov} with $\alpha=\E(V)$.
\end{theorem}
\begin{proof}
The SDE~\eqref{EmpSDE} and Lemma~\ref{lemb} show that any possible limit $(\Lambda(t))$ of the sequence $(\Lambda_N(t))$ satisfies the relation
\begin{multline}\label{eqL}
\croc{\Lambda(t),f}=\croc{\Lambda(0),f}-\int_0^t \croc{\Lambda(u),\odot f'(\cdot)}\,\diff u \\
+ \E(V)\int_0^t \croc{\Lambda(u),f'}\croc{\Lambda(u),b}\,\diff u -\int_0^t \croc{\Lambda(u),fb}\,\diff u.
\end{multline}
Proposition~\ref{pro:ExistenceUniquenessWeak} of the Appendix gives that such a process $(\Lambda(t))$ is unique, in particular it is deterministic.  The convergence in distribution of $(\Lambda_N(t))$ holds.

The process $(X_1^N(t))$ is the solution $(X(t))$ of the SDE
\begin{multline}
 X(t)= X(0)-\int_0^t X(u)\diff u\,\diff u +\Phi^N(t)-\frac{\E(V)}{N} \int_0^t X(u)b(X(u)\,\diff u
\\-\int_0^t X(u-)\int_{\R_+^2} \ind{0\leq v\leq b(X(u-))} {\cal N}_1(\diff v,\diff z,\diff u),
\end{multline}
where
\[
\Phi^N(t)=\E(V) \int_0^t \croc{\Lambda^N(u),b}\,\diff u+M_I^N(t),
\]
where $(M_1^N(t))$ is the martingale~\eqref{SDE-mart}.  The process $(\Phi^N(t))$ is converging in distribution for the uniform norm to
\[
\Phi(t)=\E(V) \int_0^t \croc{\Lambda(u),b}\,\diff u,
\]
and  the martingale vanishes with a usual argument, it is then not difficult to get that $(X_1^{N}(t))$ converges in distribution in the Skorohod's space ${\cal D}(\R_+,\R_+)$ to the solution $(\overline{X}(t))$ of the SDE
\[
 \overline{X}(t)= \left(-\overline{X}(t)+\Phi'(t)\right)\,\diff t-\int_0^t \overline{X}(t-)\int_{\R_+^2} \ind{0\leq v\leq b( \overline{X}(t-))} {\cal N}_1(\diff v,\diff z,\diff t).
\]
By using the fact that, by exchangeability, $\E(b(X_1^N(t)))=\E(\croc{\Lambda_N(t),b})$ and that $\Lambda(t)$ is deterministic, one gets that
\[
\Phi(t)=\E(\Phi(t))=\E(V) \int_0^t \E\left(\overline{X}(u)\right)\,\diff u,
\]
one concludes the proof of the theorem by using the uniqueness result of Theorem~\ref{McVTheo}. 
\end{proof}

%% file: Invariant.tex
\section{Analysis of Invariant distributions}\label{sec:Invariant}
This section is devoted to the analysis of the invariant distributions of the McKean-Vlasov process~\eqref{n2}. We shall denote in this section $\pi$ a distribution on $\R_+$ which is invariant along the McKean-Vlasov evolution. In particular the map $t\mapsto \E(b(X(t))$ is constant. If we denote $\alpha=\E(V)\E(b(X(0)))$,  the stationary process $(X(t))$ can then be seen as the solution $(Y(t))$ of the SDE
\begin{equation}\label{n1}
\diff Y(t)=(\alpha -Y(t))\diff t
-Y(t-)\int \ind{0\leq u\leq b(Y(t-))} {\cal N}(\diff u,\diff t,\diff z).
\end{equation}
Define
\[
\tau=\inf\{t> 0: Y(t)=0\}
\]
and $x(t)=\alpha(1-\exp(-t))$. 

The variable $\tau$ is the  instant of the first spike of the neuron. If $Y(0)=0$, before time $\tau$ the evolution of $(Y(t))$ is deterministic,  one has in fact  $Y(t)=x(t)$ for $t<\tau$. The Poisson property gives that
\[
\P(\tau\geq t)=\exp\left(-\int_0^t b(x(u))\,\diff u\right)
\])
The invariant distribution $\pi$ of $(Y(t))$  can then be expressed as
\[
\pi(f)=\frac{1}{\E(\tau)}\E\left(\int_0^\tau f(x(u))\,\diff u\right),
\]
for $f$ a continuous function with compact support on $[0,\alpha)$. By Fubini's Theorem
\begin{align*}
\E\left(\int_0^\tau f(x(u))\,\diff u\right)&=
\int_0^{+\infty} f(x(u))\P(\tau\geq u)\,\diff u\\
&= \int_0^{+\infty} f(x(u))\exp\left(-\int_0^u b(x(v))\,\diff v\right)\diff u\\
&= \int_0^{\alpha} \frac{f(u)}{\alpha -u}\exp\left(-\int_0^{u} \frac{b(v)}{\alpha-v}\,\diff v\right)\diff u.
\end{align*}
The measure has a compact support $[0,\alpha)$, it has finite mass if and only if
\[
\E(\tau)=\int_0^{\alpha} \frac{1}{\alpha -u}\exp\left(-\int_0^{u} \frac{b(v)}{\alpha-v}\,\diff v\right)\diff u<+\infty.
\]

\begin{theorem}\label{ThInv}
The invariant distribution of the solution of SDE~\eqref{n2} has density 
\begin{equation}\label{InvDens}
u\mapsto \frac{1}{C(\beta)(\beta \E(V) -u)}\exp\left(-\int_0^{u} \frac{b(v)}{\beta \E(V)-v}\,\diff v\right)
\end{equation}
on $[0,\beta \E(V))$, where
\[
C(\beta)=\int_0^{\beta \E(V)} \frac{1}{\beta \E(V) -u}\exp\left(-\int_0^{u} \frac{b(v)}{\beta \E(V)-v}\,\diff v\right)\diff u
\]
and $\beta$ is the solution  of the equation
\begin{equation}\label{eq:FixedPointMFE}
\beta C(\beta)=1- \exp\left(-\int_0^{\beta \E(V)} \frac{b(v)}{\beta \E(V)-v}\,\diff v\right).
\end{equation}
\end{theorem}

In the case where $b(\beta \E(V))>0$, the term 
\[
\exp\left(-\int_0^{\beta \E(V)} \frac{b(v)}{\beta \E(V)-v}\,\diff v\right)
\]
vanishes, the fixed point equation reduces to $\beta C(\beta)=1$. 

Moreover, the change of variable $x=u/\beta\E(V)$ and $y=v/\beta\E(V)$ yields the simplified formulation of $C(\beta)$:
\[C(\beta)=\int_0^{1} \frac{1}{1 -x}\exp\left(-\int_0^{x} \frac{b(\beta\E(V) y)}{1-y}\,\diff y\right)\diff x.\]	

We shall now analyze the behavior of the map $\beta\mapsto\beta C(\beta)$ in order to characterize the number of possible stationary distributions. 

\begin{lemma}\label{pro:LimitsPsi}
The function $\beta\to \beta C(\beta)$ is converging to infinity as $\beta$ gets large,  and  if 
\[
\lim_{x\to 0} \frac{b(x)}{x}= \lambda\in[0,+\infty]
\]
then
\[
\lim_{\beta\to 0} \beta C(\beta)=\frac{1}{\lambda \W}.
\]
\end{lemma}

\begin{proof}
Let us start by the behavior at infinity, by monotonicity of $b$, the quantity 
\[
B(x)=\int_0^x b(y)\,\diff y
\]
is upperbounded by $ x b(x)$. Let us fix $\delta\in (0,1)$, we have:
	\begin{align*}
		\beta C(\beta) &\geq \beta \int_0^{\delta} \frac 1 {1-x} \exp\left( -\frac 1 {1-\delta} \frac{B(\beta \W \delta)}{\beta\W} \right)\,\diff x\\
		&\geq -\beta \log(1-\delta) \exp\left(-\frac {\delta} {1-\delta} b(\beta \W \delta) \right)
	\end{align*}
	and simply taking, for arbitrary $C>0$, $\delta = C/\beta$ (for $\beta>C$), we obtain 
	\[\beta C(\beta)\geq -\beta\log(1-\frac C {\beta})\exp\left(-\frac {C} {\beta-C} b(C \W) \right).\]
	It is then easy to see that from this formula that:
	\[\lim_{\beta\to\infty}\beta C(\beta)\geq C\]
	and since $C$ is arbitrary, this precisely means that $\beta C(\beta)\to \infty$.
	
Let us now analyze the behavior of $\beta C(\beta)$ at $\beta=0$ as a function of the limit $\lambda$.
	\begin{itemize}
\item $\lambda=0$.\\ For any $\delta>0$, there exists $\beta(\delta)$ such that $b(\beta \W)/\beta \W \leq \delta$ for all $\beta \leq \beta(\delta)$, and therefore for such $\beta\leq \beta(\delta)$ we have:
		\begin{align*}
			\beta C(\beta)\geq \beta \int_0^1 \frac{1}{1-x} \exp\left(-\alpha \beta \W \int_0^x \frac{y}{1-y}\,\diff y\right)\,\diff x \geq \frac{1}{\delta \W}
		\end{align*}
		which proves that $\beta C(\beta)\to \infty$ at $\beta=0$.
\item $\lambda=\infty$.\\ For any $\delta>0$ we have for $\beta$ small enough $b(\beta \W y) \geq \delta \beta\W y$, and therefore:
		\begin{align*}
			\beta C(\beta) &\leq \beta \int_0^1 (1-x)^{\delta \W \beta -1} \exp(\delta \beta \W x)\,\diff x\\
			&= \frac 1 {\delta \W} + \beta \int_0^1 (1-x)^{\delta \W \beta}e^{\delta \beta \W x}\,\diff x
		\end{align*}
		and therefore we have $\lim_{\beta\to 0}\beta C(\beta)\leq {1}/{\delta \W}$ for arbitrarily large $\delta$, showing that $\beta C(\beta)\to 0$ at $\beta=0$.
		\item $0<\lambda<+\infty$.\\ It is easy to show using the same estimates as in the two previous cases that for any $\delta>0$ small enough we have
		\[\frac{1}{(\lambda +\delta)\W} \leq \lim_{\beta\to 0} \beta C(\beta)\leq \frac{1}{(\lambda-\delta)\W},\]
	 	which ends the proof.
	\end{itemize}
\end{proof}

From this lemma, one gets directly the following proposition.
\begin{proposition}[Number of Stationary Solutions of Mean-Field Equations]\label{cor:FixedPointsGeneral}
Let
\[
\lim_{x\to 0} \frac{b(x)}{x}= \lambda\in[0,+\infty]
\]
then, for the McKean-Vlasov process~\eqref{n2},
	\begin{itemize}
		\item if $\lambda \in [1/\W,\infty]$, then there always exists at least one non-trivial stationary 
invariant distribution. 
		\item If $\lambda \in (0,1/\W)$, the Dirac mass at $0$ is the unique  invariant distribution. 
		\item If $\lambda=0$, then the Dirac mass at $0$ is always an invariant distribution. Non-trivial invariant distribution exist if the minimal value of $\beta\to \beta C(\beta)$ is smaller than $1$. 
	\end{itemize}	
\end{proposition}
The mean-field equations can therefore present several stationary distributions. The question of the \emph{stability} of these solutions (in a sense to be made more precise) is then natural. We have seen in theorem~\ref{Th-dead} that in finite-sized networks, the only stationary solution is the trivial state in which no neuron spikes, and that the network converges towards this solution. With additional solutions arising in the mean-field limit, a deep question concern the emergence of new stationary and their significance if they are attractive. 

\subsection{Stability of the trivial solution}
In the cases where $b(0)=0$, the Dirac mass at $0$  is invariant for the mean-field equations. We investigate here the stability of this invariant distribution as a function of the local behavior of $b(x)$ at zero. The main result of the section is the following:

\begin{proposition}\label{pro:stabilityZero}
Let
	\[
\lambda =  \lim_{x\to 0} \frac {b(x)}{x} \in [0,\infty] \quad\text{ and } \rho=\lambda\E(V),\]
then
	\begin{itemize}
		\item if $\rho \in [0,1)$, the trivial solution $\delta_0$ is almost surely exponentially stable.\\ More precisely, there exists $\delta>0$ and $A_{\delta}>0$ sufficiently small such that for any initial condition with support included in $[0,A_{\delta}]$, 
		\[\limsup_{t\to \infty} \frac{\log(X_t)}{t}<-\delta \qquad a.s.\]
		\item if $\rho \in (1,\infty]$, the trivial solution is unstable in probability.\\
That is, there exists $A>0$ such that for any initial condition $X_0$ with support included an interval $[0,A]$ and mean $\mu_0>0$, there exists a $t_0>0$ such that
		\[\P\left(\sup_{t\in [0,t_0]} X_t > A\right)>0.\]
	\end{itemize}
\end{proposition}

The quantity $\rho$ can be seen as the excitation rate of the network. 
The result for $\rho<1$ is quite strong: almost any trajectory converge exponentially fast towards $0$ provided that the initial condition is chosen sufficiently close from $\delta_0$, in the sense that its support is included in a small interval around $0$. 
	
However, the instability result part of the proposition for $\rho>1$ is weaker:  whatever the initial condition, the probability of reaching in finite time a specified level away from $0$ is strictly positive. From a pathwise viewpoint, the result is indeed less strong than the exponential stability result. But from the distribution viewpoint, this corresponds to an instability of the distribution $\delta_0$ in the sense of Khasminskii~\cite{khasminskii:80}.

\begin{proof}
	Let us first deal with the case $\lambda\W<1$. Since jumps are all negative, any solution of the McKean-Vlasov equation has the upperbound:
	\[X_t\leq X_0+\int_0^t (-X_s+\W \E(b(X_s)_)\,\diff s.\]
	Moreover, if $\lambda\W<1$, there exists $\delta>0$ and $A_{\delta}>0$ such that for any $x<A_{\delta}$, 
	\[-x+\W b(x)\leq -\delta x\]
	We introduce $X_t^*=\textrm{ess sup} X_t = \inf\{u>0; \P[X_t>u]=0\}$, and assume that $X_0^*<A_{\delta}$. We show that along the evolution, the essential support of $(X_t)$ never exceeds $A_{\delta}$ and actually shrinks to ${0}$, ensuring stability of the solution $\delta_0$. Indeed, using Gronwall's lemma and the monotonicity of the map $b$, we have:
	\[X_t\leq X_0 e^{-t}+\W\int_0^t e^{t-s} \E(b(X_s))\,\diff s \leq X_0^* e^{-t}+\W\int_0^t e^{t-s} b(X_s^*)\,\diff s \]
	readily implying that:
	\[X_t^* \leq X_0 e^{-t}+\W\int_0^t e^{-(t-s)} b(X_s^*)\,\diff s.\]
	Let us now introduce the deterministic time:
	\[\tau=\inf\{t>0 \; ;\; X_t^* > A_{\delta}\}\]
	On the interval $[0,\tau)$, we have:
	\[X_t^* \leq X_0 e^{-t}+(1-\delta)\int_0^t e^{-(t-s)} X_s^*\,\diff s\] 
	and by Gronwall's lemma again, we obtain for any $t\in [0,\tau)$:
	\[X_t^* \leq X_0^* e^{-\delta t}\leq A_{\delta}e^{-\delta t}.\]
	This implies that (i) $\tau=\infty$ and (ii) $X_t^*$ converges exponentially fast towards $0$. We have therefore proved that for any initial condition with support sufficiently concentrated around $0$, the process converges almost surely towards $0$ when $t\to\infty$, hence the solution $\delta_0$ is stable. 
	
	If $\lambda\W>1$, we can find $\delta>0$ and $A_{\delta}>0$ sufficiently small so that:
	\[-x+\W b(x)-x b(x) \geq \delta x.\]
	Let $X_t$ be the solution of the McKean-Vlasov equation with initial condition $X_0$ such that $X_0^* < A_{\delta}$ (we recall that $X_t^*$ denotes in this proof the essential supremum of $X_t$). Denoting $\mu_t=\E(X_t)$, we have:
	\begin{align*}
		\mu_t &= \mu_0 + \int_0^t -\mu_s + \W \E(b(X_s))-\E(X_s\,b(X_s)) \,\diff s\\
		&=\mu_0 + \int_0^t \int_{\R} (-x+\W b(x)-x\,b(x)) \;p_s(\diff x)\,\diff s
	\end{align*}
	with $p_s$ is the distribution of $X_s$. Similarly to the previous case, let us denote $\tau$ the deterministic time:
	\[\tau = \inf\{t>0\;;\; X_t^* <A_{\delta}\}.\]
	On the interval $[0,\tau)$, we have:
	\[\mu_t \geq \mu_0 + \delta \int_0^t \mu_s\,\diff s\]
	i.e. $\mu_t\geq \mu_0 \, e^{\delta t}$. This implies that necessarily 
	\[t \leq \frac 1 \delta \log\left(\frac {A_{\delta}}{\mu_0}\right)\stackrel{\text{def.}}{=} t_0.\]
	We therefore conclude that the essential supremum of the solution exceeds $A_{\delta}$ whatever the initial condition, which means that 
	\[\P\left(\sup_{t\in [0,t_0]} X_{s} > A_{\delta}\right)>0.\]
\end{proof}

\subsection{Power firing functions}\label{sec:Power}
We now provide some specific examples, for power functions of the form $b(x)=\lambda x^{\alpha}+\delta$. We distinguish the affine ($\alpha=1$), superlinear ($\alpha>1$) and sublinear ($\alpha<1$) cases. 

\subsubsection{Affine firing functions}
We start by considering affine firing functions, and apply Proposition~\ref{cor:FixedPointsGeneral} and the characterization of the stability of Proposition~\ref{pro:stabilityZero} to study the number of invariant distributions and their stability:
\begin{proposition}[Linear Firing-Rate]\label{pro:Linear}
If $b(x)=\lambda\,x +\delta$ with $\lambda>0$ and $\rho=\lambda\W$, then if $\delta=0$
	\begin{itemize}
		\item For $\rho<1$, $\delta_0$ is the unique stationary solution and it is almost surely exponentially stable 
		\item For $\rho>1$, $\delta_0$ is unstable in probability, and there exists an additional solution to the mean-field equations. 
	\end{itemize}
For $\delta > 0$, there exists a unique, non trivial, invariant distribution.
\end{proposition}

\begin{proof}
	Affine firing functions allow analytical calculations for all quantities. Basic algebra yields
	\[C(\beta)= \frac{1}{\rho \beta +\delta} + \frac{\rho \beta }{\rho \beta+\delta}\int_0^1 (1-x)^{\rho \beta+\delta} e^{\rho \beta\,x}\,\diff x.\]
	We change variables and define $x=-1+y-\log(y)$. The map $y\mapsto -1+y-\log(y)$ is strictly decreasing on $(0,1)$ and its inverse is $\phi(x)=\exp(-W(-e^{-1-x})-1-x)$ where $W$ is the first real branch of the Lambert function $W$, see Corless \etal~\cite{corless:96}. One gets that 
	\[\beta C(\beta)=\frac{\beta}{\rho\beta+\delta} \left[1-\E\left(\Phi\left(\frac {E_1}{\rho \beta}\right)\right)\right]\]
	where $\Phi(x)=\phi^{\delta}(x)\phi'(x)$. It is then easy to see that:
	\[\frac{d}{d\beta} \beta C(\beta)=\frac{\delta}{(\rho\beta+\delta)^2}\left[1-\E\left(\Phi\left(\frac {E_1}{\rho \beta}\right)\right)\right]+\frac 1 {\rho \beta +\delta}\E\left(\frac {E_1}{\rho \beta}\Phi'\left(\frac {E_1}{\rho \beta}\right)\right).\]
	The map $\phi$ satisfies:
	\[\begin{cases}
		\displaystyle{\phi'(x)=\frac{W(-e^{-x-1})}{1+W(-e^{-x-1})}}\\
		\\
		\displaystyle{\phi''(x)=-\frac{W(-e^{-x-1})}{(1+W(-e^{-x-1}))^3}},
	\end{cases}\] 
	and therefore,
	\[\Phi'(x)=\frac{(-W(-e^{-x-1}))^{\delta+1} (\delta (1+W(-e^{-x-1}))+1)}{(1+W(-e^{-x-1}))^3}.\]
	For $x\geq 0$, $-e^{-1}\leq -e^{-x-1}\leq 0$ and hence $W(-e^{-x-1})\leq 0$, ensuring that $\varphi'(x)<0$, $\varphi''(x)>0$ and eventually $\Phi'(x)>0$. All these estimates put together prove that $\beta\mapsto \beta C(\beta)$ is strictly increasing, and therefore that there exists a unique non-trivial solution to the fixed point equation $\beta C(\beta)=1$ when $\rho>1$ or when $\delta>0$.
\end{proof}

\begin{figure}
	\centering
		\subfigure[Averaged activity]{\includegraphics[width=.4\textwidth]{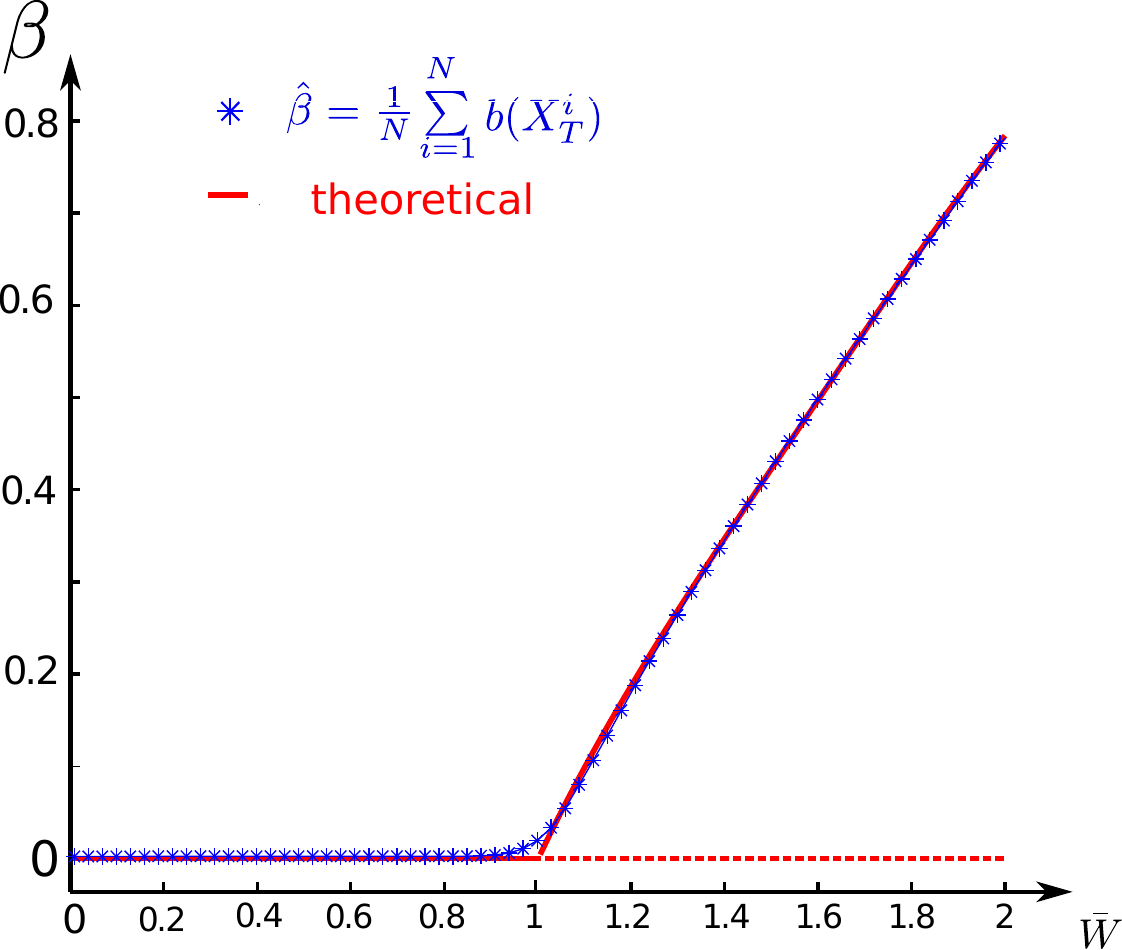}}
		\subfigure[Extinction time]{\includegraphics[width=.4\textwidth]{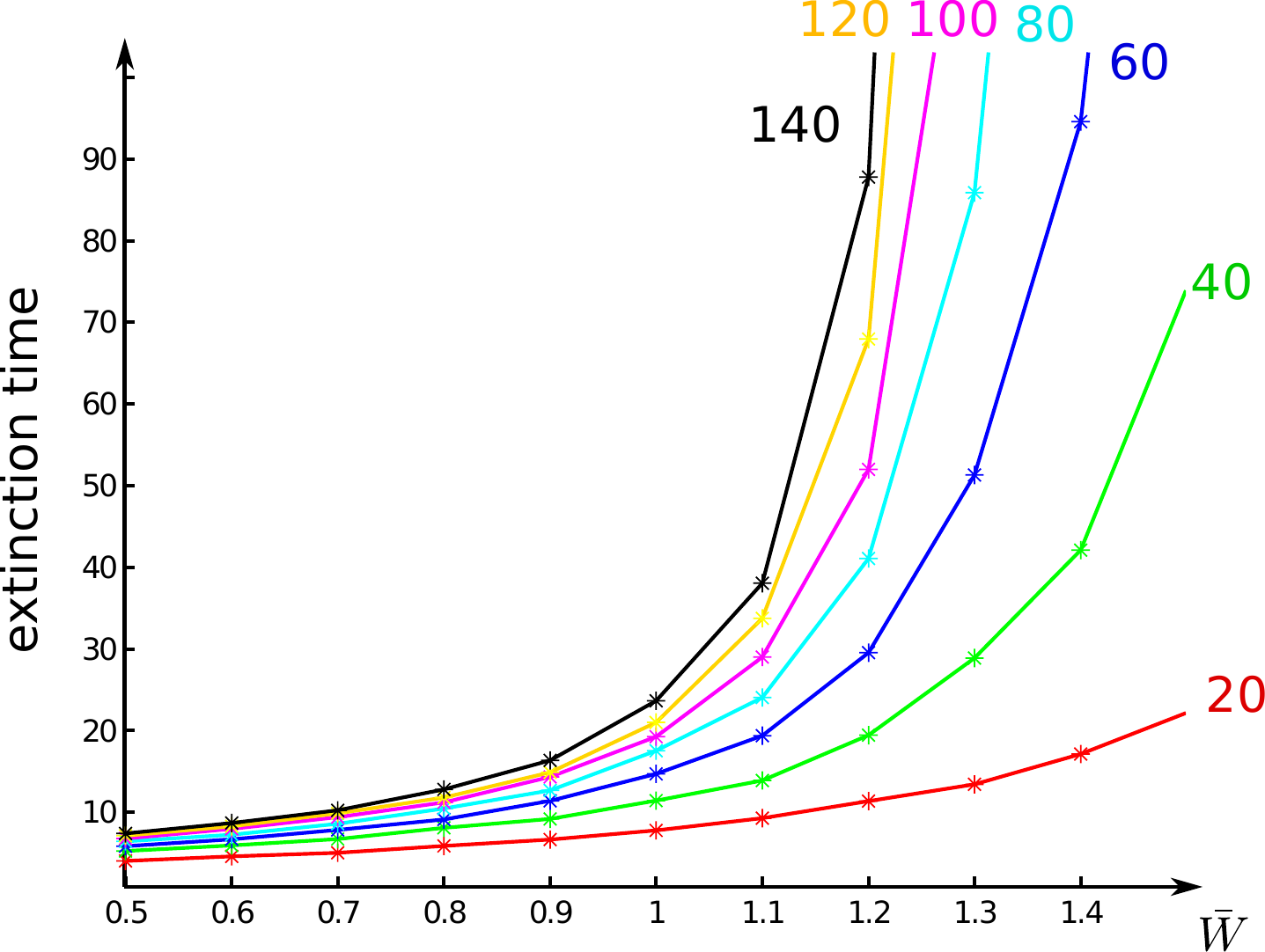}\label{fig:ExtinctionTime}}
	\caption{Affine firing-rate functions. (a) Empirical averaged activity $\hat{\beta}$ of a $2\,000$ neurons network for different values of $\W$, at time $T=100$, averaged across $30$ initial conditions uniformly drawn in $[0,1]$, with $\lambda=1$. The value of $\beta$ corresponding to the non-trivial invariant measure of the mean-field equation was computed numerically and plotted in red. (b) Extinction Time of the network as a function of $\W$ for different values of the network size $N$. }
	\label{fig:Affine}
\end{figure}

We therefore conclude that in the case $b(0)>0$, there exists a unique stationary solution, which is non-trivial, as was also the case in the finite-sized networks. 

In the case $b(0)=0$, we have shown that the only stationary distribution of finite-sized networks is the trivial solution $\delta_0$. In the mean-field limit, this solution persists whatever the value of the parameters. However, we showed that for $\rho>1$, this solution is no more stable, and non-trivial solution appears when $\rho>1$. This is what we observe in the simulations of the network (see Fig.~\ref{fig:Affine}) in the linear firing rate case $b(x)=x$ for varying values of $\W$: for $\W>1$, the trivial solution no more attracts the network, and a new solution with a non-zero value of the firing rate emerges. The value of $\beta$ at this equilibrium can be computed numerically, and shows a very good agreement with the simulations of the finite-sized network, even if this finite-sized network will eventually extinct. 

This phenomenon suggests the presence of a phase transition in the system. For small coupling ($\rho<1$), both finite-sized networks and their mean-field limit have a trivial stationary solution. In that case, the time of extinction remains small and do not dramatically depend on the network size. However, for $\rho>1$, the trivial solution is no more stable for the mean-field limit and, in that limit, a sustained activity appears. The time of extinction shows a dramatic dependence on the network size. A pseudo-stationary solution emerges, which is meta-stable for any finite system in the sense that even though the system will eventually stop firing, the time during which the system supports this non-trivial stationary firing rate diverges as the network size increases (see Fig.~\ref{fig:ExtinctionTime}). See the end of Appendix  on the simulation algorithm used to obtain the figures.

From the biological viewpoint, the non-trivial solution found corresponds to a self-sustained activity. In this regime, neurons fire independently as a Poisson processes with a common intensity. This regime is a natural regime of activity of large neuronal networks. It is a typical regime of the awake brain often referred to as the asynchronous irregular state, see Brunel~\cite{brunel:00}. 

\subsubsection{Sub-linear power firing functions}
\begin{proposition}\label{pro:sublin}
If  $b(x)=\lambda x^{\alpha}+\gamma$ with $0<\alpha<1$ and $\gamma\geq 0$, then there exists a unique non-trivial invariant distribution to the McKean-Vlasov equations. For $\gamma=0$, the trivial solution is unstable.
\end{proposition}
\begin{proof}
	For $\alpha<1$, proposition~\ref{pro:LimitsPsi} implies that  the map $\beta \mapsto \beta C(\beta)$ tends to $0$ at $\beta=0$ and to infinity when $\beta \to \infty$, ensuring the existence of a non-trivial invariant distribution. Moreover, proposition~\ref{pro:stabilityZero} shows that the trivial solution is unstable. The only result that remains to be proved is the uniqueness of the invariant distribution. To this end, we show that the map $\beta \C(\beta)$ is strictly increasing. This is done by rewriting the expression of $\beta C(\beta)$ noting $\rho=\lambda \beta^{\alpha}$ and using the expression:
	\[\beta C(\beta)=\frac{1}{\rho\beta^{\alpha-1}} + \beta \int_0^1 (1-x)^{\gamma}\frac{1-x^{\alpha}}{1-x}e^{-\rho\beta^{\alpha} \phi(x)}\,\diff x\]
	with 
	\[\phi(x)=\int_0^x \frac{y^{\alpha}-1}{1-y}\,\diff y -\log(1-x).\]
	This map is strictly increasing, tends to $0$ when $x\to 0$ and to $\infty$ at $x=1$. It is therefore invertible, and we denote $\varphi=\phi^{-1}$.
	Using the variable $z=\phi(x)$, we can express our equation as:
	\[\beta C(\beta)=\frac{1}{\rho\beta^{\alpha-1}} \left[1+\E\left(\psi_{\gamma}\left(\frac{E_1}{\rho \beta^{\alpha}}\right)\right)\right]\]
	where $E_1$ is an exponential random variable with parameter $1$ and
	\[\psi_{\gamma}=(1-\varphi)^{\gamma} \frac{1-\varphi^{\alpha}}{1-\varphi} \varphi'.\]
	 Hence, we have:
	\[\frac{d}{d\beta} \beta C(\beta) = (1-\alpha) \frac{1}{\rho \beta^{\alpha}} \left[1+\E\left(\psi_{\gamma}\left(\frac{E_1}{\rho \beta^{\alpha}}\right)\right)\right] -\frac{\alpha}{\rho^2 \beta^{2\alpha}}\E\left(E_1\psi_{\gamma}'\left(\frac{E_1}{\rho \beta^{\alpha}}\right)\right).\]
	The first term of this expression is clearly positive. The second term is handled by expressing the differential $\psi_{\gamma}'$ and showing that it is strictly negative. In details, we have:
	\[\psi_{\gamma}'=-\gamma (1-\varphi)^{\gamma-1}\psi_0 \varphi' + (1-\varphi)^{\gamma}\psi_0'\]
	and therefore we only need to show that $\psi_0'<0$. Straightforward calculations yield:
	\[\psi_0' = \left(\frac{\varphi'}{1-\varphi}\right)^2 (1-(\alpha \varphi^{\alpha-1}+(1-\alpha)\varphi^{\alpha})) + \frac{1-\varphi^{\alpha}}{1-\varphi}\varphi''.\]
	The first term is clearly positive, and the second term has the sign of $\varphi''$. Since we have:
	\[\begin{cases}
\displaystyle		\varphi'=\frac{1}{\phi'\circ \phi}, \qquad \qquad \quad \varphi''=-\frac{\varphi'}{\phi'^2\circ\phi}\phi''\circ \phi\\
\displaystyle		\phi'(x) = \frac{x^{\alpha}}{1-x}>0, \qquad \phi''(x) = \frac{\alpha x^{\alpha-1}}{1-x} + \frac{x^{\alpha}}{(1-x)^2}>0,
	\end{cases}\]
	we conclude that $\varphi''<0$. This ensures that $\beta\mapsto \beta C(\beta)$ is strictly increasing, and therefore there exists a unique invariant distribution for the mean-field equations.
\end{proof}

\subsubsection{Super-linear power firing functions}
The case $b(x)=\lambda x^{\alpha}$ with $\lambda>0$ and $\alpha>1$ shows a more intricate behavior. Proposition~\ref{pro:LimitsPsi} shows that the map $\beta\mapsto \beta C(\beta)$ diverges to infinity when $\beta\to 0$ or $\beta\to \infty$, which allowed to conclude that apart from the trivial invariant distribution, there either exist no other invariant distribution or generically an even number of non-trivial invariant distributions.

We analyze the dependence of the number of non-trivial invariant distributions as a function of the parameters. The following lemma investigates the fixed point equation $\beta C(\beta)=1$ as a function of $\lambda$, in fact of $\rho=\lambda\W^\alpha$. 
\begin{lemma}\label{lem:ManipsSuperlinear}
Denoting by $\rho=\lambda \W^{\alpha}$. For any $\beta>0$, there exists a unique $\rho(\beta)$ such that $\beta$ satisfies the fixed point equation $\beta C(\beta)=1$ and there exists $\beta_c\in\R_+^*$ such that 
	\[\rho_c = \min_{\beta\in\R_+} \rho(\beta)=\rho (\beta_c)>0.\]
\end{lemma}
\begin{proof}
	Simple algebraic manipulations allow to rewrite the fixed point equation as:
	\begin{equation}\label{eq:IPP}
		\beta C(\beta)= \frac{1}{\rho \beta^{\alpha-1}} + \beta \int_0^1 (1-x)^{\rho \beta^{\alpha}} \frac{1-x^{\alpha}}{1-x} \exp\left(-\rho \beta^{\alpha} \int_0^x \frac{y^{\alpha}-1}{1-y}\,\diff y\right)\,\diff x
	\end{equation}
	We define
	\[g(u)=\int_0^u \frac{v^\alpha}{1-v}\,\diff v\quad
\text{	and }\quad
	\Psi(x)= \int_0^1 \frac 1 {1-u} \exp(-x g(u))\,\diff u.\]
	Our fixed point equation simply reads
	\begin{equation}\label{eq:FPrho}
		\beta \Psi(\rho \beta^\alpha)=1.
	\end{equation}
	With this expression, it is now relatively easy to show that for any $\beta>0$ fixed, there exists a unique $\rho(\beta)$ such that equation~\eqref{eq:FPrho} is satisfied. Indeed, it is clear from the expression~\eqref{eq:FPrho} that $\rho\mapsto \beta \Psi(\rho \beta^\alpha)$ is strictly decreasing, tends to infinity at $\rho=0$ and to $0$ when $\rho\to\infty$. 
	
	Moreover, the map $\beta \mapsto \rho(\beta)$ has the following properties.
	\begin{itemize}
		\item Using equation~\eqref{eq:IPP}, we observe that $\beta^{\alpha-1} \rho(\beta)\to 1$ when $\beta \to 0$ hence $\rho(\beta)\to \infty$
		\item $\rho(\beta)\to \infty$ when $\beta \to \infty$. Indeed, using equation~\eqref{eq:IPP} and the series representation of $1/(1-y)$, we can show that 
		\[\lim_{\beta\to\infty} \left(\frac{\beta}{\rho(\beta)}\right)^{\frac 1 {\alpha+1}} = \int_0^{\infty}e^{-u^{\alpha+1}/{\alpha+1}}\,\diff u.\]
	\end{itemize}	
\end{proof}

The quantity $\rho_c$ constitute a transition point in the system, and governs the number of invariant distributions. The following proposition is a simple consequence of lemma~\ref{lem:ManipsSuperlinear}.
\begin{proposition}[Number of Stationary Solutions of Mean-Field Equations]\label{pro:SuperLinearFP}
	For $b(x)=\lambda x^{\alpha}$ with $\lambda>0$ and $\alpha>1$, then with $\rho_c$ and $\beta_c$ defined in Lemma~\ref{lem:ManipsSuperlinear},
	\begin{itemize}
		\item If $\rho < \rho_c$, there is no non-trivial invariant distribution
		\item If $\rho > \rho_c$, there exist at least two non-trivial invariant distributions with density defined by Relation~\eqref{InvDens} for $\beta\in\{\beta_-,\beta_+\}$ with $\beta_- <\beta_c<\beta_+$ which are solutions of the equation $\rho(\beta_{\pm})=\rho$.
			\item If $\rho =  \rho_c$, there exists a unique non-trivial invariant distribution corresponding to $\beta=\beta_c$.
	\end{itemize}
\end{proposition}

Numerical computations of the fixed point equation~\eqref{eq:FPrho} show that when $\rho>\rho_c$, there exists exactly two non-trivial invariant distributions. In order to prove this fact, we would need to show that the function $\beta\mapsto \rho(\beta)$ has a unique minimum on $\R_+$, i.e. that it is strictly decreasing on $[0,\beta_c]$ and increasing on $[\beta_c,\infty]$, or in other words that there exists a unique $\beta\in\R_+$ such that $\rho'(\beta)=0$. These conditions yield the implicit equation:
\[\Psi(x^*) = -2 x^* \Psi(x^*) \]
with $x^*=(\beta^*)^2 \rho(\beta^*)$. Showing analytically uniqueness of the solutions of this implicit equations is actually very complicated even for simple firing functions such as $b(x)=x^2$. Extensive numerical simulations tend to show however that this is the case.

\begin{figure}[htbp]
	\centering
		\includegraphics[width=.7\textwidth]{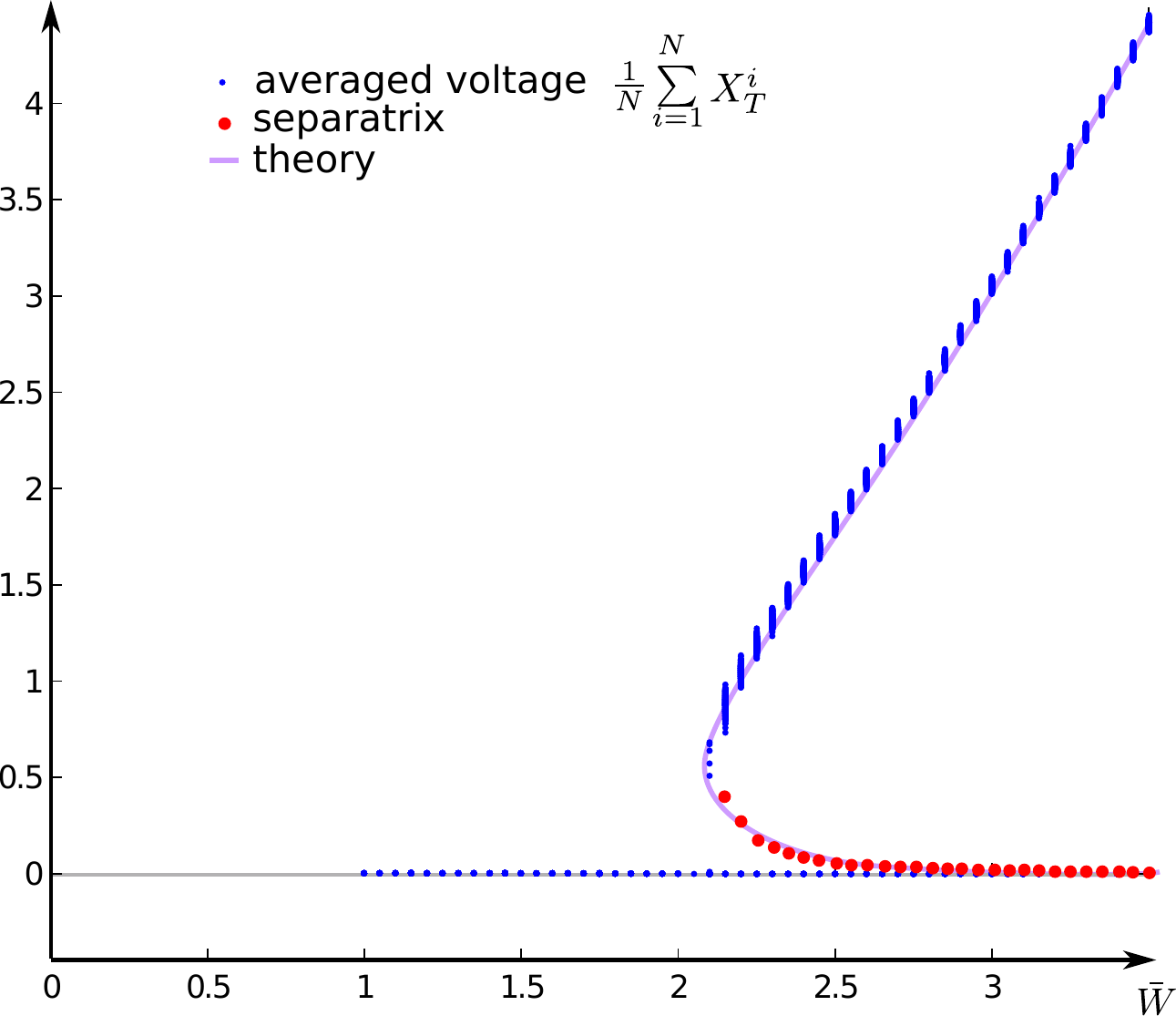}
	\caption{Quadratic firing function $b(x)=x^2$: numerical simulations of a $2\,000$ neurons network for different values of $\W$, $30$ initial conditions and $1\,000$ realizations. Each blue dot corresponds an average firing rate for a given initial condition (see text). The stationary solution $\delta_0$ persists for large values of $\W$ and the additional non-trivial invariant distribution appears when increasing $\W$. Red dots correspond to the separatrix between initial conditions converging to the trivial solution and those going to the sustained state. The purple line is the numerical solution of the mean-field equation~\eqref{eq:FixedPointMFE}, and shows a good agreement with the non-trivial and separatrix points. }
	\label{fig:Quadratic}
\end{figure}

Let us for instance discuss in more detail the case $b(x)=x^2$ (see Fig.~\ref{fig:Quadratic}). We have shown that, depending of $\W$, either the trivial distribution is the unique stationary distribution, or there exists two additional non-trivial equilibrium solutions. These solutions can be found numerically and are depicted in Fig.~\ref{fig:Quadratic} (purple line). A phase transition arises, at a specific value of $\W$, in which two additional solutions emerge. Similarly to what we did for the linear network, we extensively simulated the network in order to characterize equilibria of the system. In contrast to the linear case, all trajectories do not go to the same state, and we do expect to find certain initial conditions converging towards the trivial solution and some towards the non-trivial equilibrium. In one-dimensional dynamical systems presenting multi-stability, one typically has, between two stable equilibria, one unstable equilibrium, which acts as a separatrix, in the sense that trajectories with an initial condition on one side of the unstable equilibrium converge to the stable equilibrium on that same side. Here, the system is much more complex, and in particular it is a priori infinite-dimensional. However, one may conjecture that the limiting dynamics collapses on a smaller dimensional system. In our simulations, we considered a network made of $2\,000$ neurons, in the simple case in which the initial conditions of neurons are uniformly distributed around a value $v_0$, with a fixed standard deviation $0.2$. For $30$ fixed values of $v_0$, we simulated $1\,000$ times the network. We observed that, except for values of $\W$ close from the phase transition, that trajectories converge either towards the trivial equilibrium or towards a state with non-zero voltage. The value of $v_0$ at which a switch occurs between those trajectories going to the trivial state and those going to the non-zero state has been recorded. The average value is depicted in red in Fig~\ref{fig:Quadratic} and shows very good agreement with the middle stationary solution of the mean-field equation. The end state is characterized by one quantity per initial condition, which corresponds to the average in time (in the time interval $[90,100]$), over all neurons and over the different simulations, of the voltage variable. One point is therefore obtained for each of the $30$ initial conditions, and is depicted as a blue circle in Fig~\ref{fig:Quadratic}. The dynamics of the system, constrained to these precise initial conditions, is therefore highly similar to a one-dimensional dynamical system. Note that it is not rare that solutions of McKean-Vlasov systems reduce to low-dimensional systems. For instance in a model arising in neuroscience, it was shown in Touboul~\cite{touboulNeuralFieldsDynamics:12,touboul2012noise} that a specific, rate-based neuron model reduces exactly, in the mean-field limit, to a one-dimensional dynamical system. Although dynamics of the spiking neuron is much more complex, we conjecture that the dynamics of the firing rate is much simpler and characterized deterministically by a few statistical quantities. 

%% file: Appendix.tex
\setcounter{section}{1}
\setcounter{subsection}{0}

\section*{Appendix}
\subsection*{Some Technical Results}
\begin{lemma}\label{lemelem}
  We assume that $b$ is such that there exists $\gamma>0$ and $c>0$ such that
  \begin{equation}\label{eq:assumption2}
    b'(x) \leq \gamma b(x)+c
  \end{equation}
Then, for any $\eps>0$ and $p\in[1,3+\eps]$,  there exist a constant $\gamma_1<(3+2\eps)\gamma$,  $c_1 >0$ and a value $\eta_b>0$ such for any $a\in (0,\eta_b)$ and $x\geq 0$,
  \begin{equation}\label{eq:Discrete}
      b^p(x+a)-b^p(x) \leq a \left(\gamma_1 b^p(x)+c_1 \right).
  \end{equation}
\end{lemma}
\begin{proof}
  Let us start by noting that the inequality is trivial for $b$ bounded. We will therefore assume in the rest of the proof that $b$ diverges at infinity. We also remark that for any $1\leq p<3+\eps$, the map $b^p$ satisfies an inequality of type~\eqref{eq:assumption2} where $\gamma$ is multiplied by $p$. Indeed, for any $\delta>0$, we can find $c_{\delta}>0$ such that:
	\[
	\der{b^p(x)}{x}\leq p\gamma b^p(x) + p c b^{p-1}(x)\leq (p\gamma+\delta) b^p(x)+c_{\delta}.
	\]
	We will therefore demonstrate without loss of generality the proposition for $p=1$, and i.e. control the modulus of continuity of $b$ under assumption~\eqref{eq:assumption2}.  For an arbitrary $x_0>0$ and any $x\geq x_0$, we have:
  \begin{align*}
    \frac{b(x+a)}{b(x)} = \exp\left(\int_{x}^{x+a} \frac{b'(y)}{b(y)}\,\diff y\right)\leq e^{a \tilde{\gamma}},
  \end{align*} 
  with $\tilde{\gamma}=\gamma + {c}/{b(x_0)}$. We conclude that for $x\geq x_0$,
  \[
b(x+a)-b(x)\leq (e^{a\tilde{\gamma}}-1)b(x).
\]
The map $a\mapsto (e^{a\tilde{\gamma}}-1)/a$ is smooth, non-decreasing and tends to $\tilde{\gamma}$ at $a=0$, which can be made arbitrarily close from $\gamma$ for sufficiently large $x_0$ (since $b$ is unbounded). Therefore, there exists $x_0>0$ and $\eta>0$ such that for any $x\geq x_0$ and $a\in [0,\eta]$,
\[b(x+a)-b(x)\leq a (1+\varepsilon)\gamma b(x).\]

Denoting $c_1$ the Lipschitz constant of $b$ over the interval $[0,x_0+\eta]$, we readily obtain~\eqref{eq:Discrete} with $\gamma_1=\gamma (1+\eps)$.	 
\end{proof}
Another elementary property that is useful in our developments is the following:
\begin{proposition}\label{elmGrow}
If $x(t)$ is a non-negative c\`adl\`ag function on $\R_+$ and $\kappa>\delta>0$ such that, for $A$, $B\in\R$, 
	\[
	x(t)\leq B+ x(s)-\int_s^t x(u)^\kappa\,\diff u+ A\int_s^t x(u)^\delta\,\diff u,
	\]
holds for any $0\leq s\leq t$, then $x(t)$ is uniformly bounded on any bounded time intervals. 

Moreover, if $(x(t))$ is  $C^1$ on $\R_+$ and $B=0$, we have a uniform bound for all times:
\[
\sup_{t\geq 0} x(t)\leq C_0<+\infty,
\]
where $C_0=x(0)\wedge A^{\kappa-\delta}$.
\end{proposition}

The proof is elementary once noted that the map $x\mapsto -x^\kappa+ A x^\delta$ is upperbounded by a finite value $M>0$ and is strictly negative for any $x>A^{\kappa-\delta}$. The upperbound readily implies that $x(t)\leq x_0+B+ M \,t$. For $x$ continuously differentiable and $B=0$, the negativity of the integrand for $x>A^{\kappa-\delta}$ ensures that no trajectory exceeds $C_0=x(0) \wedge A^{\kappa-\delta}$. 	

\subsection*{Poisson Processes}
The third elementary technical result used is related to the martingales associated to marked Poisson processes. 
\begin{proposition}\label{PoisMart}
If ${\cal N}$ is a Poisson process on $\R_+^3$ with intensity measure $\diff u\otimes V(\diff z)\otimes \diff t$, $f$ is a continuous function on $\R_+^3$ and $(Y(t)=(Y_1(t),Y_2(t))$ is a c\`adl\`ag adapted processes then the process $(M(t))$ defined by 
\[
\left(\int_{s=0}^t \int_{\R_+^2} f(Y(s-),z) \left[\ind{0\leq u\leq Y_1(s-)}  {\cal N}(\diff u,\diff z,\diff s)-Y_1(s)\,\diff s\,V(\diff z)  \right]\right)
\]
is a local martingale whose previsible increasing process is given by 
\[
(\croc{M}(t))=\left(\int_0^t\diff s\int_{\R_+} \diff V(\diff z) f(Y(s),z)^2\,Y_1(s)\right).
\]
\end{proposition}
\noindent
See Rogers and Williams~\cite{rogers} and Appendix~B of Robert~\cite{Robert} for example. 

\subsection*{Uniqueness}
In the main text, we have shown tightness of the sequence of empirical measures, ensuring that the sequence is relatively compact. Moreover, we showed that the possible limits are time-dependent measures $\Lambda(t)$ that satisfy, for all $f\in C^{1}(\R)$, Equation~\eqref{eqL} that we write here as:
\begin{multline}\label{eq:weakMF}
	\Lam{t}{f}=\Lam{0}{f}+\int_0^t\left [\rule{0mm}{4mm}\Lam{u}{-x f'(x) + E(V) \Lam{u}{b} f'(x)} \right.\\\left.
	-\Lam{u}{(f(x)-f(0))\,b(x)}\rule{0mm}{4mm}\right]\,\diff u.
\end{multline}
In the above notations, $x$ is a generic symbolic variable, which we use for simplicity of notations, with the convention $\Lam{t}{f(x)} = \Lam{t}{f}$. We know that the law of the solution of the mean-field McKean-Vlasov equation~\eqref{eqintro3} satisfies this system. We aim at showing there is a unique positive Radon measure such that the nonlinear equation~\eqref{eq:weakMF} holds. We first remark that the differential equation conserves the total mass $\Lam{t}{1}=\Lam{0}{1}$. We are therefore searching for $\Lambda$ a probability measure satisfying the nonlinear equation~\eqref{eq:weakMF}. The proof of uniqueness uses the following properties:
\begin{lemma}\label{lem:BoundedSupportLamb}
	For any initial probability measure $\Lambda(0)$ of $\R_+$ with bounded support,  if  $\Lambda(t)$ is a solution of  Equation~\eqref{eq:weakMF}, there exists $C$ and $K$  such that
	\renewcommand{\theenumi}{\roman{enumi}}
	\begin{enumerate}
		\item $\displaystyle \sup_{t\geq 0} \Lam{t}{b}\leq C < \infty$.
		\item $\Lambda(t)$ has its support in  $[0,K]$ for all $t\geq 0$. 
	\end{enumerate}
\end{lemma} 
\begin{proof}
	The proof of (i) is similar to the analogous property shown on the possible solutions of the McKean-Vlasov equation. Denoting $B(t)=\Lam{t}{b}$ and using the inequality $b'(x)< \gamma b(x)+c$, we have:
	\begin{align*}
		B(t) &\leq B(0) + \int_0^t \Lam{u}{E(V) B(u) b'(x)-b(x)(b(x)-b(0))}\diff u\\
		&\leq B(0) + \int_0^t \gamma E(V) B(u)^2 +(c+b(0)) B(u) -\Lam{u}{b(x)^2}\diff u\\
		&\leq B(0) + \int_0^t (\gamma E(V)-1)B(u)^2 + (c+b(0))B(u)\diff u
	\end{align*}
	and we conclude using proposition~\ref{elmGrow}. 
	
	We now prove that any solution $\Lambda(t)$ to Equation~\eqref{eq:weakMF} has a uniformly bounded support. Let us assume that the support of $\Lambda(0)$ is contained in the interval $[0,K_0]$ and pick $f$ a continuously differentiable and non-decreasing function such that 
	\[\begin{cases}
		f(x)=0 & x<K\\
		f(x)>0 & x>K
	\end{cases}\] 
	with $K=\max(K_0, \;C\,E(V))$ with $C$ an upperbound of $\sup_{t\geq 0} \Lam{t}{b}$. Applying equation~\eqref{eq:weakMF} to $f$ and using the fact that $\Lam{t}{b}<C$, $fb\geq 0$ and $f'\geq 0$, we obtain the inequality:
	\begin{equation*}
		0\leq \Lam{t}{f} \leq \int_0^t \Lam{u}{(-x+E(V) C) f'(x)}\,\diff u \leq 0,
	\end{equation*}
	hence $\Lam{t}{f}=0$ for all $t\geq 0$, implying that the support of $\Lambda(t)$ is contained within the compact set $[0,K]$. 
\end{proof}
With these a priori estimates on $\Lambda$ in hand, we can now show the uniqueness of possible solutions to the mean-field equation. For two probability measures $\lambda_1$ and $\lambda_2$, we define the distance:
\[\Vert \lambda_1 - \lambda_2\Vert_{\mathcal{S}}=\sup\left\{\croc{\lambda_1-\lambda_2, \; f} :f\in \mathcal{S}\right\},\] 
with
\[\mathcal{S}=\Big\{f\in C^1(\R_+): \| f\|_{\infty} \vee \| f'\|_{\infty} \leq 1\Big\}\]
and note that the subset of functions of $\mathcal{S}$ with bounded support is dense in the set of continuous functions with bounded support. 

\begin{proposition}\label{pro:ExistenceUniquenessWeak}
Let $\Lambda(0)$ be a probability measure with bounded support, then Equation~\eqref{eq:weakMF} has a unique solution with initial condition $\Lambda(0)$.
\end{proposition}

\begin{proof}
	We show that $\Vert \Lambda_1(t)-\Lambda_2(t)\Vert_{\cal S} =0$ for all times. Indeed, for any $f\in \mathcal{S}$, we have, denoting $\Delta(t)=\Lambda_1(t)-\Lambda_2(t)$, 
	\[\Delt{t}{f} {=} \int_0^t \Delt{u}{{-}xf'(x){-}(f(x){-}f(0))b(x) {+}\croc{\Lambda_1(u){,}f'}b {+}\croc{\Lambda_2(u){,}b}f'}\,\diff u\]
	and therefore using the fact that $\Delta$ has a support included in the compact $[0,K]$, we have:
	\begin{equation}
		\vert \Delt{t}{f}\vert \leq \Gamma(f,K)\int_0^t \Vert \Delta(u)\Vert_{\mathcal{S}}\,\diff u
	\end{equation}
	with 
	\begin{align*}
		\Gamma(f,K)&=\sup_{x\leq K, u\geq 0} \Big \vert{-}x f'(x) {-}(f(x){-}f(0))b(x) {+} \croc{\Lambda_1(u){,}f'}b(x) {+} \croc{\Lambda_2(u){,}b}f'(x)\Big\vert\\
		& \leq \mathcal{K}:= K +C + 3 b(K)
	\end{align*}
	We therefore have for all $f\in \mathcal{S}$ the inequality:
	\[\vert \Delt{t}{f}\vert \leq \mathcal{K} \int_0^t \Vert \Delta(u)\Vert_{\mathcal{S}}\,\diff u,\]
	which is therefore also valid for the norm of $\Delta$:
	\[\Vert \Delta(t)\Vert_{\mathcal{S}}\leq \mathcal{K} \int_0^t \Vert \Delta(u)\Vert_{\mathcal{S}}\,\diff u.\]
	We conclude, by immediate recursion, that:
	\[\sup_{s\leq t}\Vert \Delta(s)\Vert_{\mathcal{S}}\leq \mathcal{K} \int_0^t \sup_{s\leq u}\Vert \Delta(s)\Vert_{\mathcal{S}}\,\diff u \leq \frac{\mathcal{K}^nt^n}{n!} \sup_{s\leq t} \Vert \Delta(s)\Vert_{\mathcal{S}},\]
	hence $\Vert \Delta(t)\Vert_{\mathcal{S}}=0$ for all $t\geq 0$. 
\end{proof}

%% file: AppendixSimu.tex
\subsection*{Simulation Algorithms}
This appendix describes the simulation algorithms used to obtain our plots for linear and quadratic rate functions in section~\ref{sec:Invariant}. We used two distinct algorithms: an exact simulation algorithm for the simulation of the extinction time, and for the sake of computational efficiency an approximate algorithm for large networks. 

The algorithm we used in order to perform efficient simulations for large networks implements the evolution of the process at discrete times $t_k=k \delta t$ with $\delta t$ a small time step. In each time interval, we compute the probability that a spike occurs within the interval. We then draw a Bernoulli random variable with this probability, and update the network state accordingly. 

This approximate dynamics allows to perform fast simulations and therefore to reach very large network size. However, computing the extinction time of the network is much more delicate. To this end, we performed, for small network sizes, exact simulations of the jump process in the case of the linear firing function $b(x)=\lambda \,x$. In the specific model we treat here, the particularly simple form of the dynamics of the variables $X_i(t)$ between spikes and the simplicity of the firing map $b$ allows to derive the cumulative density function of the spikes:
\[\P(\tau_i \geq t) = \exp(-\lambda X_i (1-\exp(-t)))\]
provided that $X_i(0)=X_i$. From this expression, one obtains the probability that neuron $i$ stops firing $p_i=\exp(-\lambda X_i)$, and also the probability of firing at time $t$ provided that the neuron does not stop firing. Therefore, although we deal with state-dependent Poisson processes, these formulae allow to simulate exactly the process and the extinction time, reached when all neurons stop firing.

%% file: main.bbl
\providecommand{\bysame}{\leavevmode\hbox to3em{\hrulefill}\thinspace}
\providecommand{\MR}{\relax\ifhmode\unskip\space\fi MR }
\providecommand{\MRhref}[2]{%
  \href{http://www.ams.org/mathscinet-getitem?mr=#1}{#2}
}
\providecommand{\href}[2]{#2}
\begin{thebibliography}{10}

\bibitem{Asmussen}
S{\o}ren Asmussen, \emph{Applied probability and queues}, John Wiley \& Sons
  Ltd., Chichester, 1987.

\bibitem{brunel:00}
N.~Brunel, \emph{Dynamics of sparsely connected networks of excitatory and
  inhibitory spiking neurons}, Journal of Computational Neuroscience \textbf{8}
  (2000), 183--208.

\bibitem{brunel:00b}
Nicolas Brunel, \emph{Dynamics of networks of randomly connected excitatory and
  inhibitory spiking neurons}, Journal of Physiology-Paris \textbf{94} (2000),
  no.~5--6, 445 -- 463.

\bibitem{burkitt:06}
AN~Burkitt, \emph{A review of the integrate-and-fire neuron model: I.
  homogeneous synaptic input}, Biological cybernetics \textbf{95} (2006),
  no.~1, 1--19.

\bibitem{burkitt:06a}
\bysame, \emph{A review of the integrate-and-fire neuron model: Ii.
  inhomogeneous synaptic input and network properties}, Biological cybernetics
  \textbf{95} (2006), no.~2, 97--112.

\bibitem{caceres-carrillo:11}
Maria C\`aceres, Jos\'e Carrillo, and Beno\^{i}t Perthame, \emph{Analysis of
  nonlinear noisy integrate \& fire neuron models: blow-up and steady states},
  The Journal of Mathematical Neuroscience (JMN) \textbf{1} (2011), no.~1
  (English).

\bibitem{caceres2014beyond}
Maria~J Caceres and Beno{\^\i}t Perthame, \emph{Beyond blow-up in excitatory
  integrate and fire neuronal networks: refractory period and spontaneous
  activity}, Journal of theoretical biology \textbf{350} (2014), 81--89.

\bibitem{chichilnisky2001simple}
EJ~Chichilnisky, \emph{A simple white noise analysis of neuronal light
  responses}, Network: Computation in Neural Systems \textbf{12} (2001), no.~2,
  199--213.

\bibitem{corless:96}
R.~M. Corless, G.~H. Gonnet, D.~E.~G. Hare, D.~J. Jeffrey, and D.~E. Knuth,
  \emph{On the {L}ambert {$W$} function}, Advances in Computational Mathematics
  \textbf{5} (1996), no.~4, 329--359. \MR{1414285 (98j:33015)}

\bibitem{Dawson:16}
Donald~A. Dawson, \emph{Measure-valued {M}arkov processes}, \'Ecole d'\'Et\'e
  de Probabilit\'es de Saint-Flour XXI---1991, Lecture Notes in Math., vol.
  1541, Springer, Berlin, 1993, pp.~1--260.

\bibitem{demasi-presutti:14}
A.~De~Masi, A.~Galves, E.~L{\"o}cherbach, and E.~Presutti, \emph{Hydrodynamic
  limit for interacting neurons}, Arxiv preprint arXiv:1401.4264, February
  2014.

\bibitem{tanre:12}
F.~Delarue, J.~Inglis, S.~Rubenthaler, and E.~Tanr{\'e}, \emph{Global
  solvability of a networked integrate-and-fire model of mckean-vlasov type},
  Annals of Applied Probability (2015), To Appear.

\bibitem{delarue:hal-01001716}
Fran{\c c}ois Delarue, James Inglis, Sylvain Rubenthaler, and Etienne
  Tanr{\'e}, \emph{{Particle systems with a singular mean-field
  self-excitation. Application to neuronal networks.}}, {Stochastic Processes
  and Applications} (2015), 40.

\bibitem{gerstein-mandelbrot:64}
Georges~L. Gerstein and Benoit Mandelbrot, \emph{Random walk models for the
  spike activity of a single neuron}, Biophysical Journal \textbf{4} (1964),
  41--68.

\bibitem{khasminskii:80}
RZ~{Has'minskii}, \emph{Stochastic stability of differential equations}, Kluwer
  Academic Pub, 1980.

\bibitem{inglis2014mean}
James Inglis and Denis Talay, \emph{Mean-field limit of a stochastic particle
  system smoothly interacting through threshold hitting-times and applications
  to neural networks with dendritic component}, arXiv preprint arXiv:1409.8221
  (2014).

\bibitem{knight:72}
B.~W. Knight, \emph{Dynamics of encoding in a population of neurons}, J. Gen.
  Physiol. \textbf{59} (1972), 734--766.

\bibitem{lapicque:07}
L~Lapicque, \emph{Recherches quantitatifs sur l'excitation des nerfs traitee
  comme une polarisation}, J. Physiol. Paris \textbf{9} (1907), 620--635.

\bibitem{Levin}
David~A. Levin, Yuval Peres, and Elizabeth~L. Wilmer, \emph{Markov chains and
  mixing times}, American Mathematical Society, Providence, RI, 2009.

\bibitem{Loynes}
R.M. Loynes, \emph{The stability of queues with non independent inter-arrival
  and service times}, Proc. Cambridge Ph. Soc. \textbf{58} (1962), 497--520.

\bibitem{FL2014}
Fournier N. and E.~L{\"o}cherbach, \emph{On a toy model of interacting
  neurons}, Arxiv preprint arXiv:1410.3263, October 2014.

\bibitem{Nummelin}
Esa Nummelin, \emph{General irreducible {M}arkov chains and nonnegative
  operators}, Cambridge University Press, Cambridge, 1984.

\bibitem{pakdaman2010dynamics}
Khashayar Pakdaman, Beno{\^\i}t Perthame, and Delphine Salort, \emph{Dynamics
  of a structured neuron population}, Nonlinearity \textbf{23} (2010), no.~1,
  55.

\bibitem{pakdaman2013relaxation}
\bysame, \emph{Relaxation and self-sustained oscillations in the time elapsed
  neuron network model}, SIAM Journal on Applied Mathematics \textbf{73}
  (2013), no.~3, 1260--1279.

\bibitem{pakdaman2012adaptation}
Khashayar Pakdaman, Beno{\^\i}t Perthame, Delphine Salort, et~al.,
  \emph{Adaptation and fatigue model for neuron networks and large time
  asymptotics in a nonlinear fragmentation equation},  (2012).

\bibitem{pillow2005prediction}
Jonathan~W Pillow, Liam Paninski, Valerie~J Uzzell, Eero~P Simoncelli, and
  EJ~Chichilnisky, \emph{Prediction and decoding of retinal ganglion cell
  responses with a probabilistic spiking model}, The Journal of Neuroscience
  \textbf{25} (2005), no.~47, 11003--11013.

\bibitem{pillow2008spatio}
Jonathan~W Pillow, Jonathon Shlens, Liam Paninski, Alexander Sher, Alan~M
  Litke, EJ~Chichilnisky, and Eero~P Simoncelli, \emph{Spatio-temporal
  correlations and visual signalling in a complete neuronal population}, Nature
  \textbf{454} (2008), no.~7207, 995--999.

\bibitem{Robert}
Philippe Robert, \emph{Stochastic networks and queues}, Stochastic Modelling
  and Applied Probability Series, vol.~52, Springer, New-York, June 2003.

\bibitem{rogers}
L.~C.~G. Rogers and David Williams, \emph{Diffusions, {M}arkov processes, and
  martingales. {V}ol. 2: It\^{o} calculus}, John Wiley \& Sons Inc., New York,
  1987.

\bibitem{rolls-deco:10}
ET~Rolls and G~Deco, \emph{The noisy brain: stochastic dynamics as a principle
  of brain function}, Oxford university press, 2010.

\bibitem{Scheutzow}
M.~Scheutzow, \emph{Periodic behavior of the stochastic brusselator in the
  mean-field limit}, Probability Theory and Related Fields \textbf{72} (1986),
  425--462.

\bibitem{stein:65}
R.~B. Stein, \emph{A theoretical analysis of neuronal variability}, Biophysics
  Journal \textbf{5} (1965), 173--194.

\bibitem{Sznitman}
A.S. Sznitman, \emph{Topics in propagation of chaos}, {\'E}cole d'{\'E}t\'{e}
  de {P}robabilit\'{e}s de {S}aint-{F}lour XIX --- 1989, Lecture Notes in
  Maths, vol. 1464, Springer-Verlag, 1991, pp.~167--243.

\bibitem{touboulNeuralFieldsDynamics:12}
Jonathan Touboul, \emph{Mean-field equations for stochastic firing-rate neural
  fields with delays: derivation and noise-induced transitions}, Physica D:
  Nonlinear Phenomena \textbf{241} (2012), no.~15, 1223---1244.

\bibitem{touboulneuralfields:14}
Jonathan Touboul, \emph{The propagation of chaos in neural fields}, Annals of
  Applied Probability \textbf{24} (2014), no.~3, 1298--1328.

\bibitem{touboulSingular:14}
\bysame, \emph{Spatially extended networks with singular multi-scale
  connectivity patterns}, Journal of Statistical Physics \textbf{156} (2014),
  no.~3, 546--573 (English).

\bibitem{touboul2012noise}
Jonathan Touboul, Geoffroy Hermann, and Olivier Faugeras, \emph{Noise-induced
  behaviors in neural mean field dynamics}, SIAM Journal on Applied Dynamical
  Systems \textbf{11} (2012), no.~1, 49--81.

\end{thebibliography}
